\documentclass[11pt]{article}

\usepackage{authblk}
\usepackage{graphicx}
\usepackage[utf8]{inputenc}
\usepackage[numbers]{natbib}
\usepackage{amsmath}
\usepackage{amsfonts}
\usepackage{amssymb}
\usepackage{caption}
\usepackage{subcaption}
\usepackage[spaces,hyphens]{url}

\title{Controlling Fractional Difference Equations Using Feedback.}
\author{Divya D. Joshi\\ Department of Physics, Rashtrasant Tukadoji Maharaj Nagpur University, Amravati Rd, Nagpur-
440033, Maharashtra, India
	\and Sachin Bhalekar\\ School of Mathematics and Statics, University of Hyderabad, Lingampally CR Rao Road, Hyderabad Central University Rd, Gachibowli, Hyderabad-500046, Telangana, India
	\and Prashant M. Gade\\ Department of Physics, Rashtrasant Tukadoji Maharaj Nagpur University, Amravati Rd, Nagpur-
	440033, Maharashtra, India}

\date{\today}
\newtheorem{The}{Theorem}[section]

\newtheorem{Lem}{Lemma}[section]
\newtheorem{Ex}{Example}[section]
\newtheorem{Def}{Definition}[section]
\newtheorem{Pro}{Property}[section]

\begin{document}
\maketitle

% Optional TOC
% \tableofcontents
% \pagebreak

%--Paper--

\begin{abstract}
	One of the most popular methods of controlling dynamical systems is feedback. It can be used without acquiring detailed knowledge of the underlying system. In this work, we study the stability of fractional-order linear difference equations under feedback. The stability results are derived for an arbitrary feedback time $\tau$. We study the cases of $\tau=1$ and $\tau=2$ in further detail. The extension to the stability of fixed points under feedback for nonlinear fractional order difference equations with fixed points $ x_{*}=0$ is also carried out. 
\end{abstract}

%\begin{keywords}
%	Fractional order map,  Delay Control, Stability
%\end{keywords}

%--/Paper--

\section{Introduction}
Differential Equations have been used for modeling various phenomena in natural sciences for a long time. This modeling helps us in understanding the physical phenomena and controlling them if necessary. The differential equations have been found useful in  modeling plethora 
of systems. They range from the spreading of diseases \cite{smith2004sir}, emotional self-regulation in romantic couples \cite{steele2011latent}, economic development and growth \cite{zhang2005differential} to tumor growth \cite{enderling2014mathematical}. This modeling is inadequate in certain systems and generalization is required. Fractional differential equations are the generalization of differential equations  for systems with memory. The order of fractional differential equations could be real or even complex. Even though fractional calculus has been around for more than 300 years, fractional differential equations are applied in real world situations only in the past few decades. We find the first mention of fractional derivatives by Leibniz	and L'Hospital as early as 1695 but the very first definition of fractional derivative was introduced by Liouville and Riemann in the second half of the 19th century. Thereafter several eminent mathematicians like Caputo, Hadamard, Grunwald, Letnikov, Riesz, and others gave various definitions for fractional derivatives. This development in the field of calculus opened the door for many other fields where fractional differential equations are an essential part of mathematical modeling.

The systems defined by fractional differential equations are said to be non-local. The reason is that the future depends on its entire history. Thus, systems governed by fractional differential equations have long-term memory.  
Hence, fractional differential equations are ideal models for systems where memory plays an important role. Such systems are found in diverse fields. In the field of material science, fractional differential equations are used to define viscoelastic materials and their order determines the amount of viscosity and elasticity present in the materials \cite{schiessel1995generalized}. Fractional order epidemic models of various diseases like COVID-19 \cite{zhang2020dynamics}, Ebola \cite{dokuyucu2020fractional}, HIV \cite{naik2020global}, influenza A (H1N1) \cite{gonzalez2014fractional} have shown promising results that helped to understand, analyze and control the spread of these diseases. Seismological studies have shown that fractional order intensity measures for probabilistic seismic demand modeling applied to highway bridges have superior  performance  compared to  traditionally used intensity measures. This improved efficiency and proficiency, at the same time, maintaining practicality and sufficiency \cite{shafieezadeh2012fractional} for seismological applications. Several books have been devoted to applications of fractional calculus in real-life applications \cite{uchaikin2013fractional} and recent applications of fractional calculus in the field of science and engineering can be found in \cite{SUN2018213}. 

Fractional order systems are different from integer order systems in several respects. However, they show phenomena observed in integer order systems such as chaotic or aperiodic behavior. In certain systems, chaos is not desired and several control	schemes are designed to control the chaos. Two of the most popular control schemes are the Ott-Grebogi-Yorke method \cite{ott1990controlling} and the Pyragas method \cite{pyragas1992continuous}. For the first method, we need to know the stable and unstable manifold of the desired orbit.  For fractional order systems, the presence of such manifolds itself is a matter of debate \cite{bhalekar2020nonexistence}. The method that can be used without any detailed	knowledge of the systems is the feedback method suggested by Pyragas. We study the possibility of controlling fixed points using feedback and show that the method indeed works in the fractional case as well.

In the case of systems with delay $\tau$, the value of the current state depends on its value $\tau$ steps back. Such systems have applications both in modeling as well as control. Fractional differential equations with delay have found applications in control theory. Control theory deals with the control of dynamical systems, uses feedback in several cases and has wide applications. Controllers may be designed using a feedback mechanism. Here, the output of the system is fed back to the system through a controller to influence the behavior of the system and give the desired output. 

Stability analysis of delay differential equation with control is carried out by Bhalekar \cite{bhalekar2016stability}. Controllers have been designed for controlling the fractional ordered systems with input delay 	\cite{si2009sliding,yin2012design,you2020relative,nirmala2016controllability}. The stability of the Cournot duopoly model with distributed time delays is discussed by Culda {\it{et al.}} in \cite{culda2022stability}. The existence of chaotic behavior, stability, and its synchronization of fractional ordered systems with delay have been studied for several systems like the Ikeda system \cite{jun2006chaotic}, logistic system \cite{wang2008chaos} and Chen system \cite{daftardar2012dynamics}. The ongoing scenario in the world demands newer and updated knowledge of the spreading of diseases and their controllability. As discussed earlier, fractional-order differential equations are used for modeling epidemiological problems. Certain realistic situations demand the introduction of a delay in resulting equations and could lead to better modeling and prediction. A rich dynamical behavior is obtained  for the infection model of fractional order with delay \cite{latha2017fractional}.

The fractional differential equations with delay have been studied in various contexts. Rihan {\it{et al.}} investigated the fractional order model of the interactions of tumor cells and the immune systems with two different time delays. The stability of the solutions was observed to have improved and the model leads to various complex behaviors \cite{rihan2020dynamics}. Alzahrani {\it{et al.}} studied the adverse effect of untimely or delayed reporting of infectious diseases more realistically by using fractional differential equations with delay \cite{alzahrani2022repercussions}.  

Fractional difference equations are a relatively much less studied models. Difference equations are the discretized versions of differential equations. The study of fractional difference equations can be seen as an approach to studying fractional differential equations by the finite difference method. Furthermore, if an integer-order difference equation is generalized to a fractional-order one, then it is able to model the memory properties in the system. This is due to the non-local property of the fractional order difference operator. In this generalization, all the values of the system from an initial point are considered while evaluating the new value. These equations demand fewer computational resources and are simpler to code. The stability conditions for difference equations of fractional order and even complex order have already been obtained by \cite{stanislawski2013stability,bhalekar2022stability,joshi2022study}. Stability results of two-term fractional difference equations are proposed by Brandibur and Kaslik \cite{brandibur2022stability}. Stability conditions are crucial for studying the control and synchronization of the systems. For systems of complex order, we have numerically investigated a fractional difference equation along with a delay term which can be viewed as a controller \cite{joshi2022study}. We observed that the parameter range over which chaos is obtained is reduced on the introduction of delay term for the complex fractional Lozi map. This can be viewed as a control.  One of the simplest control is a case of a fixed point where the system gives steady output. Though the above work indicates that feedback can be useful for control in fractional order maps, we need rigorous analytic conditions for practical applications. 

The theory of fractional finite differences is initiated by Lubich \cite{lubich1986discretized} and Miller and Ross \cite{miller1988fractional}. The topic is further developed by Atici and coworkers \cite{atici2007transform}. The stability analysis of these equations is presented in \cite{vcermak2015explicit, brandibur2022stability}.
The chaotic systems of these kinds are employed in the image processing by Abdeljawad {\it{et al.}} \cite{abdeljawad2020discrete}. Atici and Sengul introduced the fractional calculus of variations and derived the discrete Euler-Lagrange equation  \cite{atici2010modeling}. A tumor-growth in cancer is modelled by using the nabla-fractional difference equations in \cite{atici2019study}. Ouannas, Batiha and Pham \cite{ouannas2023fractional} proposed the theory and applications of chaotic difference equations of fractional order. These systems are recently used to model COVID-19 \cite{djenina2022novel}. Fractional order Mandelbrot set and Julia sets are studied in \cite{danca2023mandelbrot}.

This paper gives the stability analysis for stable fixed points for systems defined by fractional difference equations coupled with a delay term. The exact analysis is carried out for the linear system and is extended to the nonlinear maps. Stabilization or destabilization of fixed points of fractional order maps with feedback is studied. Because this is a discrete system, the control term can be added proportionally to the value at the previous time. This is unlike the differential equations studied above where the control term is within the integration and thus the control is also fractional order. This is a simpler and more practical case where the control at a given time depends on the value of variable $\tau$ steps back and not the entire history. We obtain the stability bounds for this system. The analysis is essentially for linear systems. However, we observe that the same analysis gives equally good stability bounds for nonlinear maps with appropriate linearization near the zero fixed point. In this work, we have considered the $h-$ difference operator while defining the system. For simplicity, we have taken $h=1$. Therefore, $\tau$ has integer values. Fractional values of $\tau$ can be considered for $h\ne1$.

The plan of the paper is as follows. We give essential definitions followed by the model. We carry out stability analysis of the fractional order difference equations with delay. We show that the stability conditions can be expressed in an equivalent matrix form. We study the cases $\tau=1$ and $\tau=2$ in greater detail. The case of large $\tau$ is also discussed briefly. Several examples are given for linear and nonlinear systems to corroborate our results with numerical evidence in various systems. 

\section{Preliminaries}
In this section, we have given some basic definitions. Let $h>0$, $a\in \mathbb{R}$, $(h\mathbb{N})_a = \{a,a+h,a+2h,\ldots\}$ and $\mathbb{N_\circ} = \{0,1,2,\ldots\}$. 
\begin{Def}(see \cite{mozyrska2015transform})
	The Z-transform of a sequence $ \{y(n)\}_{n=0}^\infty $ is a complex function given by
	\begin{equation*}
		Y(z)=Z[y](z)=\sum_{k=0}^{\infty} y(k) z^{-k}
	\end{equation*}
	where $z \in \mathbb{C}$ is a complex number for which the series converges absolutely. 
\end{Def}
\begin{Def}(see \cite{ferreira2011fractional,bastos2011discrete})
	Let $ h > 0 ,\; a \in \mathbb{R}$ and 
	$ (h\mathbb{N})_a = \{ a, a+h, a+2h, \ldots\} $.
	For a function $x : (h\mathbb{N})_a \rightarrow  \mathbb{C}$, 
	the forward h-difference operator is defined as 
	$$(\Delta_h x)(t)=\frac{x(t+h)- x(t)}{h},$$
	where $t$	$ \in (h\mathbb{N})_a $.
\end{Def}
Throughout this article, we take $a = 0$ and 
$h = 1$. We write $\Delta$ for $\Delta_1 $.
Now, we generalize the fractional order operators 
defined in [\cite{mozyrska2015transform,ferreira2011fractional,bastos2011discrete}].
\begin{Def}
	For a function  $x : (h\mathbb{N})_a \rightarrow  \mathbb{C}$ 
	the fractional h-sum of order 
	$\alpha = u +\iota v \in \mathbb{C}, u>0 $ is given by
	\begin{equation*}
		(_{a}\Delta_h^{-\alpha}x)(t) 
		= \frac{h^\alpha}{\Gamma(\alpha)}\sum_{s=0}^{n}\frac{\Gamma(\alpha+n-s)}{\Gamma(n-s+1)} x(a+sh),\\
	\end{equation*}
	where, $t=a+(\alpha+n)h, \; n \in \mathbb{N_\circ}$.
\end{Def}
For $h=1$ and $a=0$	, we have
\begin{eqnarray*}
	(\Delta^{-\alpha}x)(t) 	&=&\frac{1}{\Gamma(\alpha)}\sum_{s=0}^{n}\frac{\Gamma(\alpha+n-s)}{\Gamma(n-s+1)}x(s)\\
	&=&\sum_{s=0}^{n}
	\left(
	\begin{array}{c}
		n-s+\alpha-1\\
		n-s\\
	\end{array}
	\right)
	x(s).
\end{eqnarray*}                         
Here, we used the generalized binomial coefficient
\begin{equation*}
	\left(
	\begin{array}{c}
		\mu \\
		\eta\\
	\end{array}
	\right)
	=\frac{\Gamma(\mu+1)}{\Gamma(\eta+1)\Gamma(\mu-\eta+1)},\\
	\; \mu ,\eta \in \mathbb{C}, \; \text{Re}(\mu)>0,\;\text{and Re}(\eta)>0.
\end{equation*}
If $n$ $\in$ $\mathbb{N_\circ}$ then
\begin{eqnarray*}
	\left(
	\begin{array}{c}
		\mu \\
		n
	\end{array}
	\right)
	=\frac{(\mu + 1)}{n!\Gamma(\mu-n+1)}
	=\frac{\mu(\mu-1)\ldots(\mu-n-1)}{n!}.
\end{eqnarray*}
\begin{Def}
	For $n \in \mathbb{N_\circ}$ and $\alpha=u+\iota v \in \mathbb{C}, u>0,$ we define
	\begin{eqnarray*}
		\tilde{\phi}_{\alpha}(n)=
		\left(
		\begin{array}{c}
			n+\alpha-1\\
			n\\
		\end{array}
		\right)
		=(-1)^n
		\left(
		\begin{array}{c}
			-\alpha\\
			n
		\end{array}
		\right).
	\end{eqnarray*}
\end{Def}
\textbf{Note}: The convolution $\tilde{\phi}_{\alpha}*x$ of the sequences $\tilde{\phi}_{\alpha}$ and $x$ is defined as
\begin{equation*}
	\left(\tilde{\phi}_{\alpha}*x\right)(n)=\sum_{s=0}^{n}\tilde{\phi}_{\alpha}(n-s)x(s)
\end{equation*}
\begin{equation*}
	\therefore (\Delta^{-\alpha}x)(n)=(\tilde{\phi}_{\alpha}*x)(n)\\.
\end{equation*}
\begin{eqnarray*}
	\therefore Z(\Delta^{-\alpha}x)(n)=Z\left(\tilde{\phi}(n)\right)Z(x(n))\\
	=(1-z^{-1})^{-\alpha}X(z),	
\end{eqnarray*}
where $X$ is $Z$ transform of $x$.
\begin{Pro}(see \cite{petalez})
	The time-shifting property shows how a change in the discrete function's time domain alters the Z-domain.
\end{Pro}
\begin{equation*}
	Z[x(k-n)]=z^{-n}X(z)
\end{equation*}
\textbf{Proof}:
From Definition 2.1 we have,
\begin{equation*}
	X(z)=\sum_{k=0}^{\infty}x(k)z^{-k}.
\end{equation*}
Consider $k-n=m$ i.e., $k=m+n$. Thus, we write the z-transform equation as 
\begin{equation*}
	Z[x(k-n)]=\sum_{k=0}^{\infty}x(k-n)z^{-k}\\
	=\sum_{m=0}^{\infty}x(m)z^{-(m+n)}\\
	=\sum_{m=0}^{\infty}x(m)z^{-m}z^{-n}\\
	=z^{-n}\sum_{m=0}^{\infty}x(m)z^{-m}\\
	=z^{-n}X(z)
\end{equation*}
\begin{Lem}
	For $\alpha \in \mathbb{C},\; \text{Re}(\alpha)>0$,
	\begin{equation*}
		Z(\tilde{\phi}_{\alpha}(t))=\frac{1}{(1-z^{-1})^{\alpha}}.
	\end{equation*}
\end{Lem}
\textbf{Proof}: We have, 
\begin{eqnarray*}
	Z(\tilde{\phi}_{\alpha}(t))&=&\sum_{j=0}^{\infty}\tilde{\phi}_{\alpha}(j)z^{-j}\\
	&=&\sum_{j=0}^{\infty}\left(
	\begin{array}{c}
		j+\alpha-1\\
		j
	\end{array}
	\right)z^{-j}\\
	&=&\sum_{j=0}^{\infty}(-1)^{j}\left(
	\begin{array}{c}
		-\alpha\\
		j
	\end{array}
	\right)z^{-j}\\
	&=&(1-z^{-1})^{-\alpha}.
\end{eqnarray*}
by using Newton's generalization of the Binomial Theorem. 
\cite{niven1969formal,link1}.
\section{Model}
Consider the fractional order linear difference equation
\begin{equation}
	x(t)=x_0+\sum_{j=0}^{t-1}\frac{\Gamma(t-j+\alpha-1)}{\Gamma(\alpha)\Gamma(t-j)}(a-1)x(j). \label{eqn1}
\end{equation}
In this paper, we study the stability analysis of a fractional difference equation (\ref{eqn1}) coupled with a delay term. 
\subsection*{Modeling Equation}
Introducing the delay term in the equation as $bx(t-\tau)$, we get
\begin{equation}
	x(t)=b x(t-\tau)+x_{0}+\sum_{j=0}^{t-1}\frac{\Gamma(t-j+\alpha-1)}{\Gamma(\alpha)\Gamma(t-j)}\left((a-1)x(j)\right), \label{eqn2}
\end{equation}
where, $b \in \mathbb{R}$ and $a \in \mathbb{C}$.
\begin{equation*}
	\therefore	x(t+1)-b x(t+1-\tau)=x_0+(a-1)(\tilde{\phi}_{\alpha}*x)(t).
\end{equation*}
Taking Z-transform on both sides, we get
\begin{eqnarray}
	zX(z)-zx_0-b z^{(1-\tau)}X(z) \nonumber	&=&\frac{x_0}{1-z^{-1}}+\frac{(a-1)}{(1-z^{-1})^{\alpha}}X(z). \nonumber\\	
	\therefore	X(z)\left(z-b z^{(1-\tau)}-\frac{(a-1)}{(1-z^{-1})^{\alpha}}\right)&=&zx_0+\frac{x_0}{1-z^{-1}},\label{eqn3}
\end{eqnarray} where $|z|<1$.
\subsection*{Characteristic Equation}
From (\ref{eqn3}), the characteristic equation of (\ref{eqn2}) is
\begin{equation}
	\left(z(1-z^{-1})^{\alpha}-b(1-z^{-1})^{\alpha}z^{(1-\tau)}-(a-1)\right)=0, \label{eqn4}
\end{equation}
where the condition $|z|<1$ should be satisfied. 
Putting $z=e^{(it)},$ we get,
\begin{eqnarray}
	e^{(it)}(1-e^{-it})^{\alpha}-b(1-e^{-it})^{\alpha}e^{(1-\tau)it}-(a-1)&=&0,  \nonumber \\
	\text{i.e. }	e^{(it)}(1-e^{-it})^{\alpha}-b(1-e^{-it})^{\alpha}e^{(1-\tau)it}+1&=& a.\label{eqn 5}
\end{eqnarray}
\subsection*{Matrix Representation}
Equation (\ref{eqn2}) can be represented equivalently as the following system
\begin{eqnarray}
	x(t)&=&x_0+\sum_{j=0}^{t-1}\frac{\Gamma(t-j+\alpha-1)}{\Gamma(\alpha)\Gamma(t-j)}\left((a-1)x(j)\right)+b y(t), y(t)=x(t-1) \text{ for } \tau=1.\label{eqn 6}\\
	\therefore x(t+1)&=&x_0+(a-1)(\tilde{\phi}_{\alpha}*x)(t)+b y(t+1). \nonumber
\end{eqnarray}
Taking z-transform, we get
\begin{equation*}
	zX(z)-zx_0	=\frac{x_0}{1-z^{-1}}+\frac{(a-1)}{(1-z^{-1})^{\alpha}}X(z)+b(zY(z)-zy_0), Y(z)=\frac{X(z)}{z}+x(-1).
\end{equation*}
\begin{equation*}
	\therefore	z(1-z^{-1})^{\alpha}X(z)-(a-1)X(z)-bz(1-z^{-1})^{\alpha}Y(z) = x_0(1-z^{-1})^{\alpha -1} +z(1-z^{-1})^{\alpha}x_0-bzy_0(1-z^{-1})^{\alpha},
\end{equation*}
\begin{eqnarray*}
	zY(z)-X(z)=zx(-1).\\
	\therefore
	\begin{bmatrix}
		z(1-z^{-1})^{\alpha}-(a-1) & -bz(1-z^{-1})^{\alpha}\\
		-1 & z
	\end{bmatrix}
	\begin{bmatrix}
		X(z)\\
		Y(z)
	\end{bmatrix}
	=0.\\
	\therefore
	\begin{vmatrix}
		z(1-z^{-1})^{\alpha}-(a-1) & -bz(1-z^{-1})^{\alpha}\\
		-1 & z
	\end{vmatrix}
	=0.\\
	\text{i.e. }	z(1-z^{-1})^{\alpha}-b(1-z^{-1})^{\alpha}+1=a.
\end{eqnarray*}
This is the equation for the model with $\tau=1$.
Similarly, for $\tau=2$ we get,
\begin{eqnarray*}
	\begin{vmatrix}
		z(1-z^{-1})^{\alpha}-(a-1) & 0 & -bz(1-z^{-1})^{\alpha}\\
		-1 & (z-1)+1 & 0\\
		0 & -1 & (z-1)+1
	\end{vmatrix}
	=0.\\
\end{eqnarray*}
We can generalize and write $\tau+1$ dimensional determinant for delay $\tau$  as follows.

\begin{eqnarray*}
	\begin{vmatrix}
		z(1-z^{-1})^{\alpha}-(a-1) & 0 & 0 & $\ldots$ & -bz(1-z^{-1})^{\alpha}\\
		-1 & (z-1)+1 & 0 & \ldots & 0 \\
		0 & -1 & (z-1)+1  & \ldots & 0 \\
		\vdots & \ddots & \ddots & \ldots & \vdots \\
		0 & \ldots & -1  & (z-1)+1 & 0\\
		0 & 0 & \ldots  & -1 & (z-1)+1
	\end{vmatrix}
	=0.\\
\end{eqnarray*}

\subsection*{Boundary curve}
Separating real and imaginary parts of equation (\ref{eqn 5}), we get $$Re(a)=2^{\alpha}\left(\sin\left(\frac{t}{2}\right)\right)^{\alpha}\left(\cos\left(\frac{\alpha\pi}{2}+t\left(1-\frac{\alpha}{2}\right)\right)-b\cos\left(\frac{\alpha\pi}{2}+t\left(1-\tau-\frac{\alpha}{2}\right) \right)\right)+1,$$ $$Im(a)=2^{\alpha}\left(\sin\left(\frac{t}{2}\right)\right)^{\alpha}\left(\sin\left(\frac{\alpha\pi}{2}+t\left(1-\frac{\alpha}{2}\right)\right)-b\sin\left(\frac{\alpha\pi}{2}+t\left(1-\tau-\frac{\alpha}{2}\right) \right)\right).$$
Therefore, the parametric representation of the boundary condition is
\begin{equation}
	\begin{aligned}
		\beta(t)= & \{2^{\alpha}\left(\sin\left(\frac{t}{2}\right)\right)^{\alpha}\left(\cos\left(\frac{\alpha\pi}{2}+t\left(1-\frac{\alpha}{2}\right)\right)-b\cos\left(\frac{\alpha\pi}{2}+t\left(1-\tau-\frac{\alpha}{2}\right) \right)\right)+1, \\
		& 2^{\alpha}\left(\sin\left(\frac{t}{2}\right)\right)^{\alpha}\left(\sin\left(\frac{\alpha\pi}{2}+t\left(1-\frac{\alpha}{2}\right)\right)-b\sin\left(\frac{\alpha\pi}{2}+t\left(1-\tau-\frac{\alpha}{2}\right) \right)\right)\} ,
	\end{aligned}\label{eqn 7}
\end{equation}
for $t \in [0,2\pi]$ in the complex plane.
If the complex number $a$ lies inside this anticlockwise oriented simple closed curve $\beta(t)$ then the system (\ref{eqn2}) will be asymptotically stable.
\subsection{Stability Analysis with $\tau=1$}
\begin{figure*}[h]
	\centering
	\includegraphics[scale=0.25]{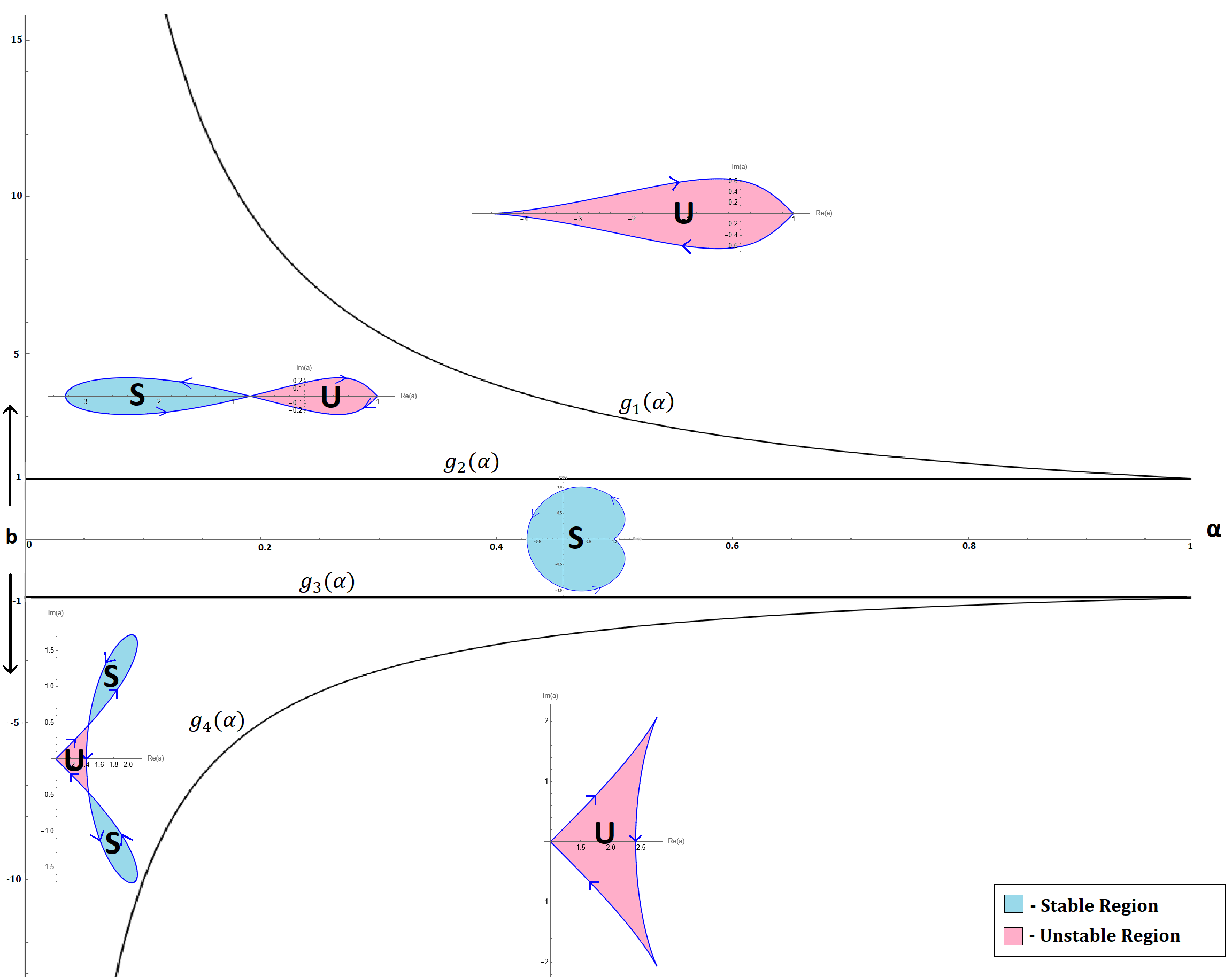}
	\caption{bifurcation regions in $b \alpha$-plane for 
		$\tau=1$ with representative stability diagrams. Here we have magnified the range on the y-axis for representation purposes. (NOT TO SCALE) }
	\label{fig1}
\end{figure*}
\begin{The}
	Consider the delayed fractional order equation (\ref{eqn2}) with $\tau=1$. Let $b=g_{1}(\alpha), b=g_{2}(\alpha)$ and $b=g_{4}(\alpha)$ be the bifurcation curves as shown in Figure \ref{fig1} in $b \alpha$-plane which are the branches of the implicit curve
	\begin{eqnarray}
		g(b,\alpha)&=& b\alpha \cos\left(\frac{1}{2}(1+\alpha)\left(\pi-\arctan\left(\frac{R}{S}\right)\right)\right)\nonumber\\
		&&-(-1+\alpha)\cos\left(\frac{1}{2}\left(\pi(1+\alpha)-(-1+\alpha)\arctan\left(\frac{R}{S}\right)\right)\right)\nonumber\\
		&&+\sin\left(\frac{1}{2}\left(\pi \alpha-(-3+\alpha)\arctan\left(\frac{R}{S}\right)\right)\right), \label{impl}
	\end{eqnarray}
	where,
	\begin{equation*}
		R=\frac{-1+\alpha+b\alpha-b\alpha^{2}
			-\sqrt{(1+b\alpha)^{2}(1-2\alpha+4b\alpha+\alpha^{2})}}{b\alpha},
	\end{equation*}
	and
	\begin{eqnarray*}
		S&=&-\frac{\sqrt{2}}{b\alpha}(-(1+2(-1+b)b^{2}\alpha^{3}
		+b^{2}\alpha^{4}+\alpha^{2}-3b^{2}\alpha^{2}-2\alpha+2b\alpha\\
		&&+(1-\alpha-b\alpha+b\alpha^{2})(\sqrt{(1+b\alpha^{2})
			(1-2\alpha+4b\alpha+\alpha^{2})})))^{\frac{1}{2}}.
	\end{eqnarray*}
	Furthermore, define the curve $b=g_3(\alpha)$ as the straight line $b=-1$ in the $b \alpha $-plane.
\end{The}
We have the following stability results:\\
1. If $b \in (-\infty,g_{4}(\alpha)) \cup (g_{1}(\alpha),\infty)$,
then the system (\ref{eqn2}) is unstable.\\
2. If $g_{4}(\alpha)<b<g_{3}(\alpha)$, then the boundary curve (\ref{eqn 7})
produces three disjoint regions; two of which are stable and
one is unstable.\\
3. If $g_{3}(\alpha)<b<g_{2}(\alpha)$, then the boundary curve (\ref{eqn 7})
generates a bounded region which is a stable region for the system.\\
4. If $g_{2}(\alpha)<b<g_{1}(\alpha)$, then the single stable
region in case 3 gets divided into two regions; one is stable
and the other is unstable.\\
5. If $b=g_{2}(\alpha)$, then the curve (\ref{eqn 7}) is a smooth
curve and the inside part is stable.\\
6. If $b=g_{3}(\alpha)$, then $\beta(0)=\beta(\pi)$ and the
stable region gets divided into two parts.\\
7. If $b=g_{1}(\alpha)$ or $b=g_{4}(\alpha)$, then the curve
$\beta$  will have cusps.\\
\textbf{Proof:}
The boundary curve $\beta(t)$ defined by (\ref{eqn 7}) is closed because the initial point $\beta(0)$ has the same value as that of final point $\beta(2\pi)$. If $\beta(t)$ is simple and its orientation is anticlockwise (positive) then the bounded part inside this curve is the stable region for the system (\ref{eqn2}). As the parameter values change, the simple closed curve $\beta(t)$ transforms to the non-simple curve i.e. a curve containing multiple points. Therefore, the region bounded by the branch of the non-simple curve  $\beta(t)$ with anticlockwise orientation is the required stable region for the system (\ref{eqn2}).\\

It is observed that the formation of cusps in $\beta(t)$ is responsible for the generation of multiple points (self-intersections).
Thus, the parameter values at which cusps are produced in $\beta(t)$ are the points of bifurcation.
At the cuspidal point $\beta(t_{0})$, the curve $\beta(t)=(x(t),y(t))$ cannot have the derivative $\beta'(t)=\frac{y'(t)}{x'(t)}$. Here, $x(t)$ and $y(t)$ are the coordinates of the curve $\beta(t)$. In this case, $\lim_{t\rightarrow t_{0}}\beta'(t)$ can take either $\frac{0}{0}$ or $\frac{\infty}{\infty}$ form. It is observed that at initial point $t=0$, this limit takes $\frac{\infty}{\infty}$ form. Therefore, the initial point $\beta(0)$ is responsible for the bifurcations in some cases.\\

To check another possibility (i.e. $\frac{0}{0}$ form), we solve the equations $x'(t)=0$ and $y'(t)=0$ simultaneously. We eliminate $t$ between these equations and get a relation in the parameters $b$ and $\alpha$. This is carried out by	squaring and adding the equations $x'(t)=0$ and $y'(t)=0$. Now, we get the value of $t$ as functions of $b$ and $\alpha$. Substituting this $t$ value in the equation $y'(t)=0$, we get the implicit expression (\ref{impl}).

This curve $g(b,\alpha)=0$ has 3 branches viz $b=g_{1}(\alpha), b=g_{2}(\alpha)$ and $b=g_{4}(\alpha)$  as shown in Figure \ref{fig1}. Note that,  $b=g_{2}(\alpha)$ is the straight line $b=1$ in the $b \alpha$-plane.

Furthermore, the curve $\beta(t)$ intersects itself at $\beta(0)$ and $\beta(\pi)$ as the parameter $b$ passes through the value $-1$. At this intersection point $\beta(0)=\beta(\pi)$, $\beta'(t)$ does not exist at $t=0,\pi$ and we have a bifurcation apart from cusps. Thus, $b=-1$ is also a bifurcation curve in the $b-\alpha $ plane. We denote it by $b=g_{3}(\alpha)$.

If the $b$ value is above the curve $b=g_{1}(\alpha)$, the orientation of the boundary curve $\beta(t)$ becomes clockwise and hence the system becomes	unstable for the values of '$a$' in the bounded as well as unbounded regions described by $\beta(t)$.

As $b$ passes through the curve $b=g_{1}(\alpha)$, the $\beta(t)$ becomes multiple curves and generates stable and unstable bounded regions as shown in Figure \ref{fig1}.	If we keep on decreasing the value of $b$ the unstable part becomes	smaller and vanishes at $b=1$. Note that $\beta$ is smooth at $t=0$ for $g_{1}(\alpha)<b\le g_{2}(\alpha)$. If $g_{2}(\alpha)<b<g_{3}(\alpha)$ then the $\beta(t)$ is a simple closed curve with an anticlockwise orientation. Therefore, the region bounded by $\beta(t)$ is the stable region and the unbounded part of the complex plane is the unstable region.

As we take the values of $b$ between the curves $b=g_{3}(\alpha)$ and $b=g_{4}(\alpha)$, $\beta(t)$ produces two stable and one unstable bounded regions as shown in the Figure \ref{fig1}.	The size of the unstable region increases if we decrease the value of $b$.

At the bifurcation curve $b=g_{4}(\alpha)$, the stable regions vanish completely and the $\beta(t)$ becomes a simple closed curve with clockwise orientation. Thus, the system is completely unstable for $b<g_{3}(\alpha)$. This proves the result.

\subsection{Illustrative Examples}
\begin{figure*}[h]
	\begin{subfigure}{\columnwidth}
		\centering
		\includegraphics[scale=0.55]{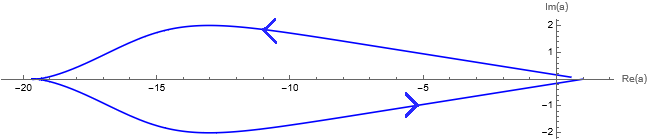}    
		\caption{Stability curve with $\tau=1$, $\alpha=0.1$ and $b=18.3$}
		\label{fig2a}
	\end{subfigure}
	\begin{subfigure}{0.45\textwidth}
		\includegraphics[width=0.7\columnwidth]{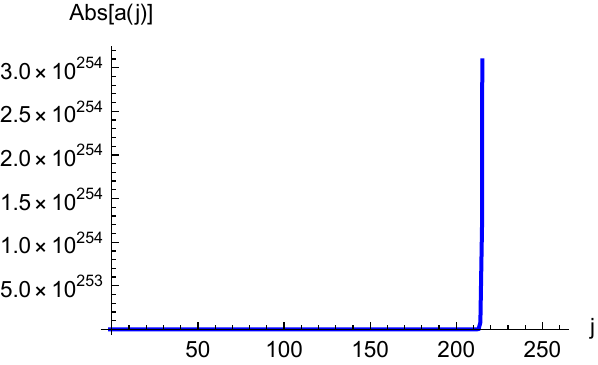}  
		\caption{Unstable solution for $a=-2$.}
		\label{fig2b}
	\end{subfigure}
	\begin{subfigure}{0.45\textwidth}
		\includegraphics[width=0.7\columnwidth]{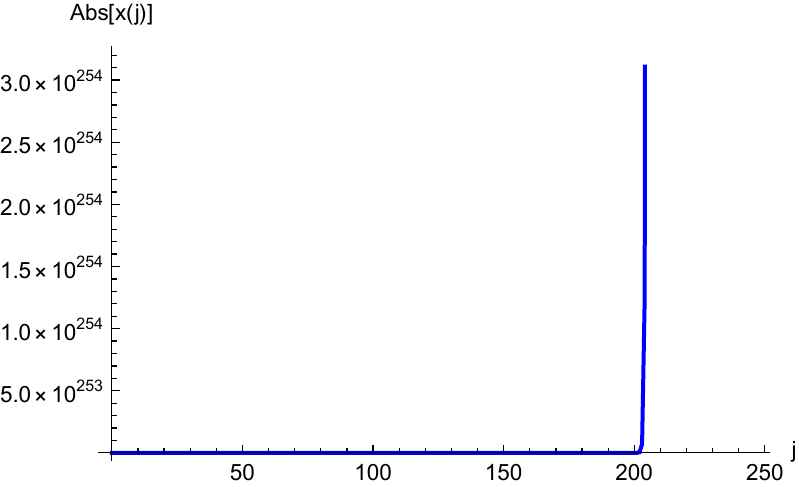}  
		\caption{Unstable solution for $a=(0.4+0.4\iota)$.}
		\label{fig2c}
	\end{subfigure}
	\caption{Example 1}
	\label{fig2}
\end{figure*}

\begin{figure*}[h]
	\begin{subfigure}{\columnwidth}
		\centering
		\includegraphics[scale=0.55]{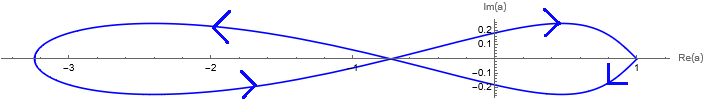}   
		\caption{Stability curve with $\tau=1$, $\alpha=0.5$ and $b=2$}
		\label{fig3a}
	\end{subfigure}
	\begin{subfigure}{0.45\textwidth}
		\includegraphics[width=0.7\columnwidth]{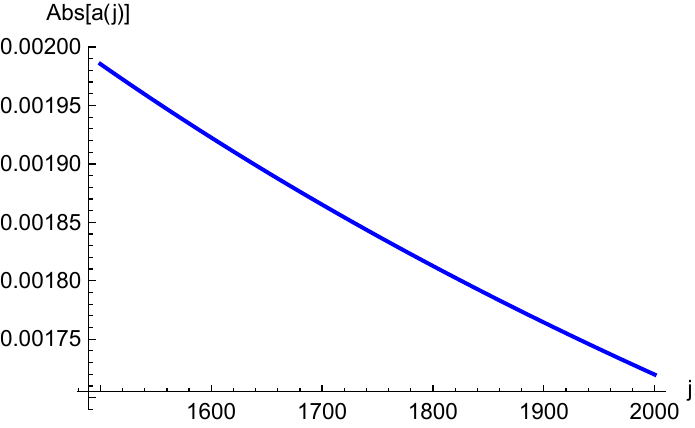}
		\caption{Stable solution for $a=-1.2$.}
		\label{fig3b}
	\end{subfigure}
	\begin{subfigure}{0.45\textwidth}
		\includegraphics[width=0.7\columnwidth]{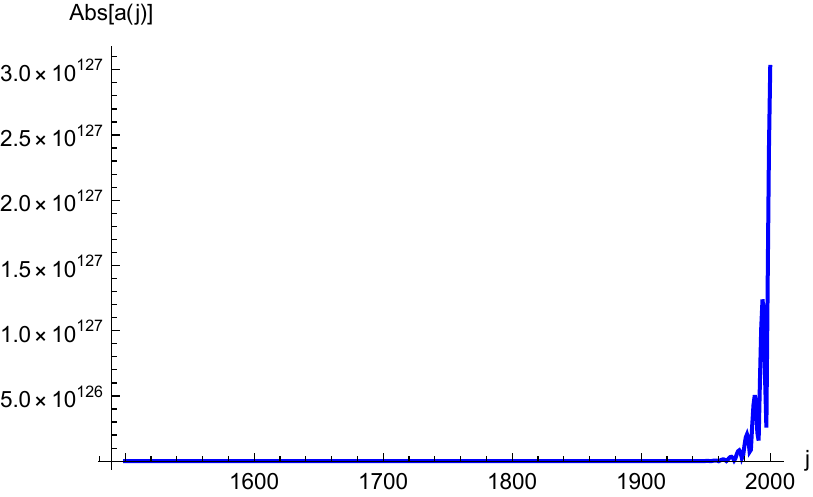}  
		\caption{Unstable solution for $a=0.2$.}
		\label{fig3c}
	\end{subfigure}
	\caption{Example 2}
	\label{fig3}
\end{figure*}

\begin{figure*}[h]
	\begin{subfigure}{\columnwidth}
		\centering
		\includegraphics[scale=0.45]{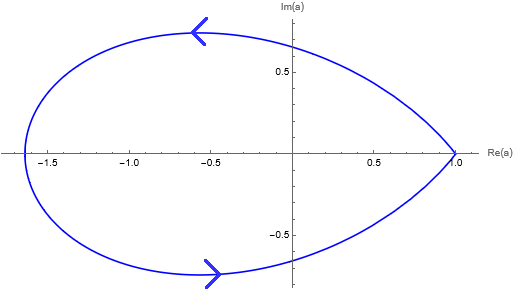} 
		\caption{Stability curve with $\tau=1$, $\alpha=0.4$ and $b=1$}
		\label{fig4a}
	\end{subfigure}
	\begin{subfigure}{0.45\textwidth}
		\includegraphics[width=0.7\columnwidth]{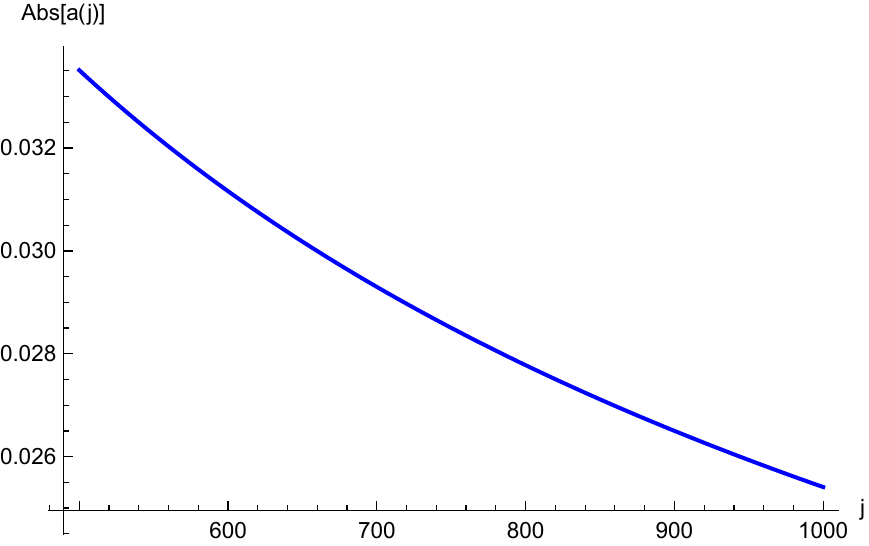}  
		\caption{Stable solution for $a=0.5$.}
		\label{fig4b}
	\end{subfigure}
	\begin{subfigure}{0.45\textwidth}
		\includegraphics[width=0.7\columnwidth]{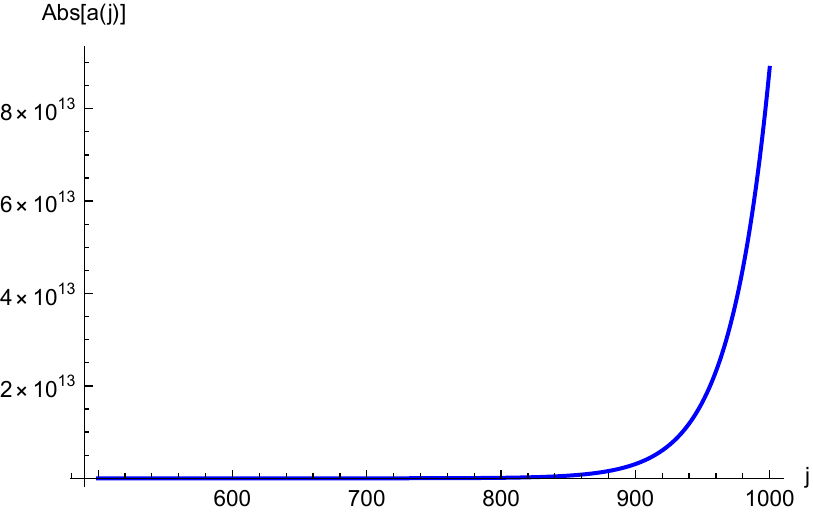}  
		\caption{Unstable solution for $a=-0.7\iota$.}
		\label{fig4c}
	\end{subfigure}
	\caption{Example 3}
	\label{fig4}
\end{figure*}

\begin{figure*}[h]
	\begin{subfigure}{\columnwidth}
		\centering
		\includegraphics[scale=0.65]{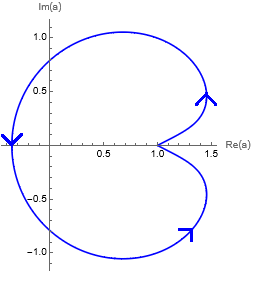} 
		\caption{Stability curve with $\tau=1$, $\alpha=0.3$ and $b=0.1$}
		\label{fig5a}
	\end{subfigure}
	\begin{subfigure}{0.45\textwidth}
		\includegraphics[width=0.65\columnwidth]{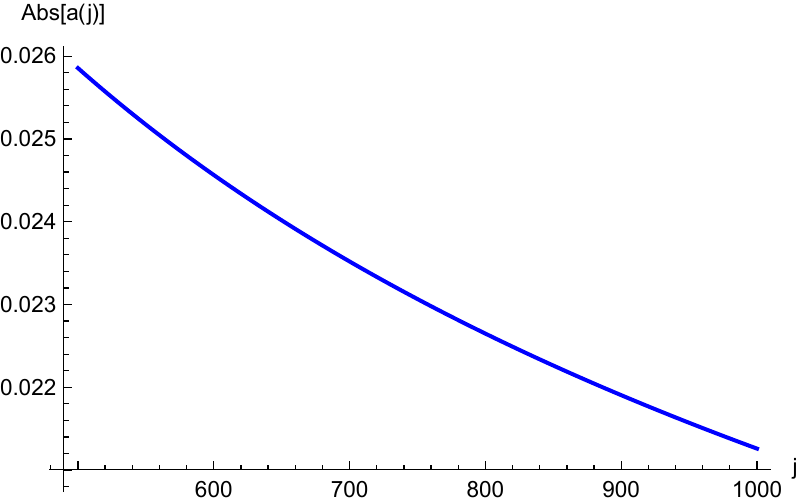}  
		\caption{Stable solution for $a=-0.3$.}
		\label{fig5b}
	\end{subfigure}
	\begin{subfigure}{0.45\textwidth}
		\includegraphics[width=0.65\columnwidth]{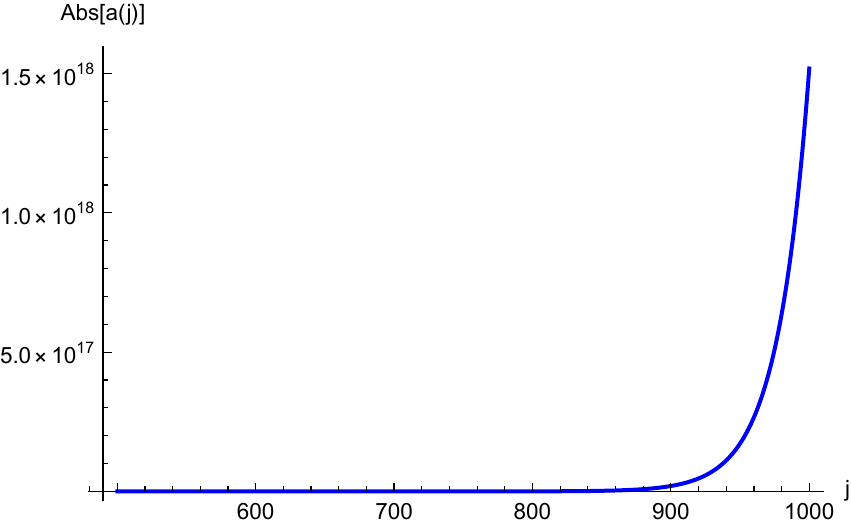}  
		\caption{Unstable solution for $a=-0.4$.}
		\label{fig5c}
	\end{subfigure}
	\caption{Example 4}
	\label{fig5}
\end{figure*}

\begin{figure*}[h]
	\begin{subfigure}{\columnwidth}
		\centering
		\includegraphics[scale=0.5]{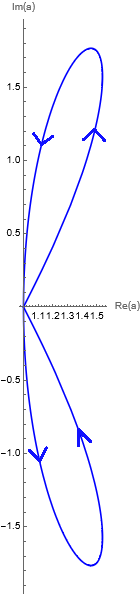}  
		\caption{Stability curve with $\tau=1$, $\alpha=0.7$ and $b=-1$}
		\label{fig6a}
	\end{subfigure}
	\begin{subfigure}{0.45\textwidth}
		\includegraphics[width=0.7\columnwidth]{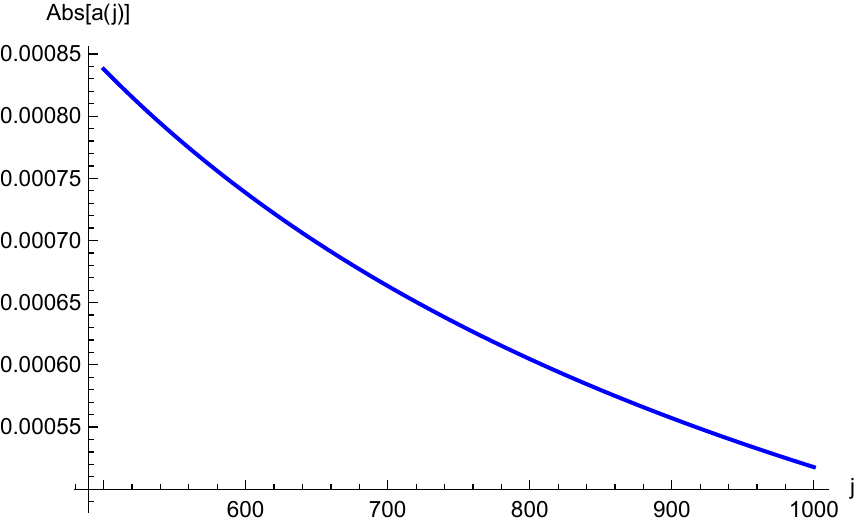}  
		\caption{Stable solution for $a=1.1+0.5\iota$.}
		\label{fig6b}
	\end{subfigure}
	\begin{subfigure}{0.45\textwidth}
		\includegraphics[width=0.7\columnwidth]{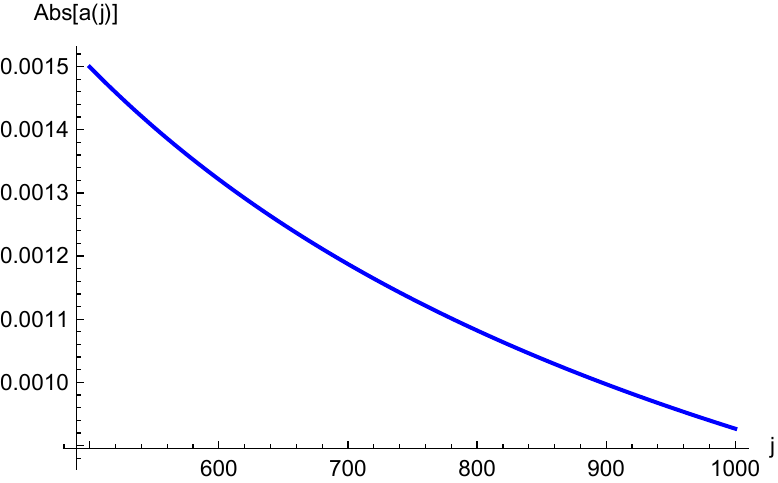}  
		\caption{Stable solution for $a=1.3-0.8\iota$.}
		\label{fig6c}
	\end{subfigure}
	\caption{Example 5}
	\label{fig6}
\end{figure*}

\begin{figure*}[h]
	\begin{subfigure}{\columnwidth}
		\centering
		\includegraphics[scale=0.65]{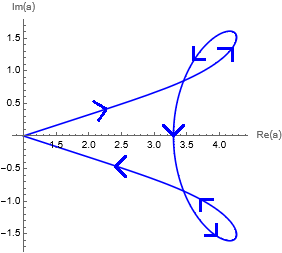}    
		\caption{Stability curve with $\tau=1$, $\alpha=0.2$ and $b=-3$}
		\label{fig7a}
	\end{subfigure}
	\begin{subfigure}{0.45\textwidth}
		\includegraphics[width=0.7\columnwidth]{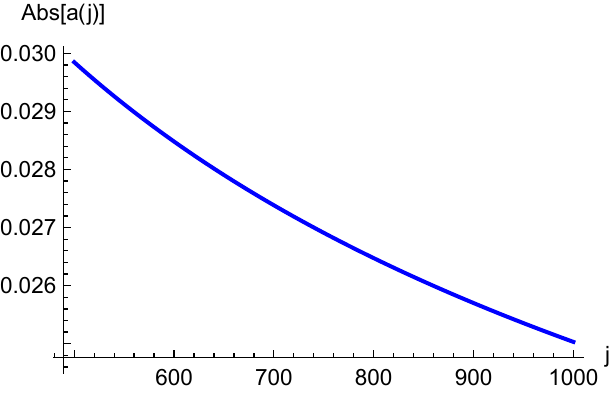}  
		\caption{Stable solution for $a=4+1.3\iota$.}
		\label{fig7b}
	\end{subfigure}
	\begin{subfigure}{0.45\textwidth}
		\includegraphics[width=0.7\columnwidth]{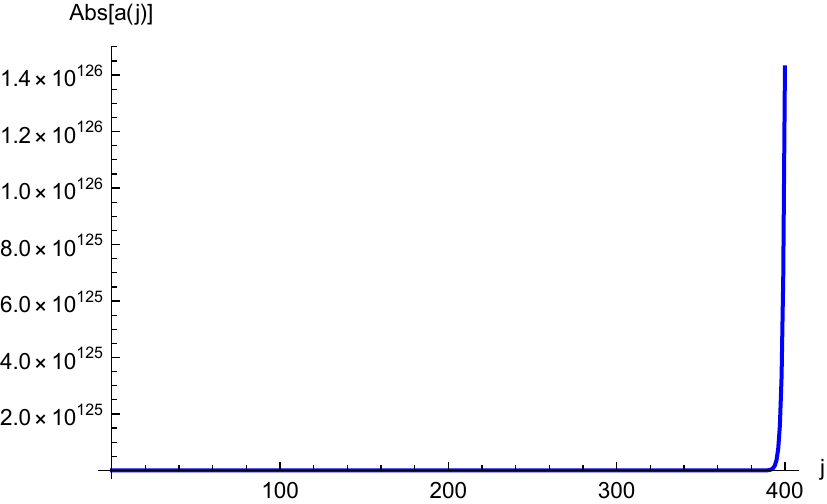}  
		\caption{Unstable solution for $a=2$.}
		\label{fig7c}
	\end{subfigure}
	\caption{Example 6}
	\label{fig7}
\end{figure*}

\begin{figure*}[h]
	\begin{subfigure}{\columnwidth}
		\centering
		\includegraphics[scale=0.57]{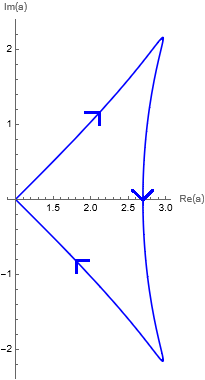} 
		\caption{Stability curve with $\tau=1$, $\alpha=0.5$ and $b=-2.2$}
		\label{fig8a}
	\end{subfigure}
	\begin{subfigure}{0.45\textwidth}
		\includegraphics[width=0.75\columnwidth]{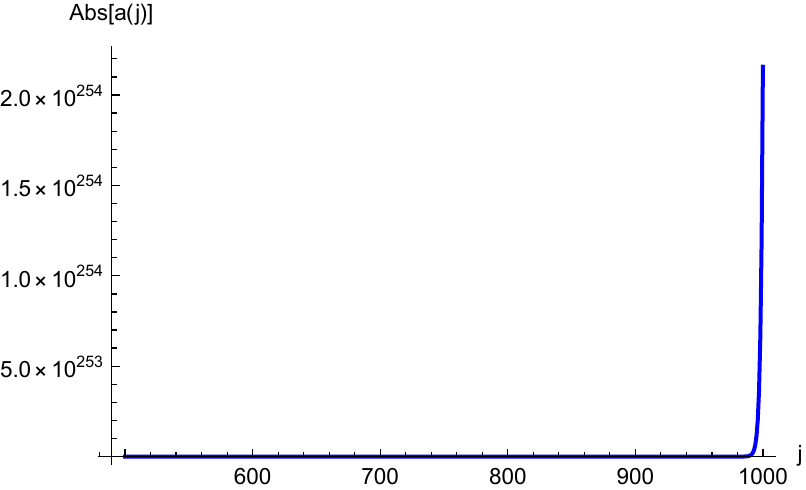}  
		\caption{Unstable solution for $a=1.5$.}
		\label{fig8b}
	\end{subfigure}
	\begin{subfigure}{0.45\textwidth}
		\includegraphics[width=0.75\columnwidth]{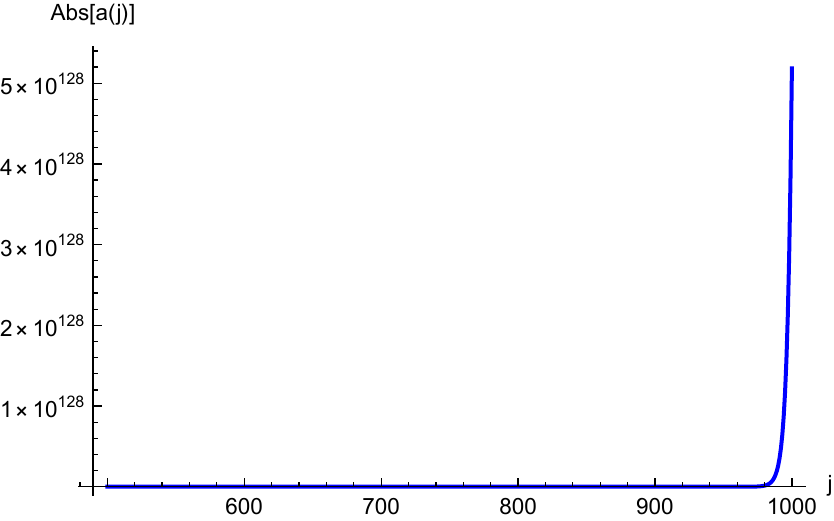}  
		\caption{Unstable solution for $a=2.8$.}
		\label{fig8c}
	\end{subfigure}
	\caption{Example 7}
	\label{fig8}
\end{figure*}

In this section, we provide corroborating evidence
for stability results obtained in the above section for $\tau=1$
and various values of $b$.
We study  various regions 
shown in Figure \ref{fig1}. 
\begin{Ex}
	Let us begin with the  region that lies above 
	$g_{1}(\alpha)$. We take $b=18.3$ and 
	$\alpha=0.1$ (Figure \ref{fig2a}). Here, we 
	observe a non-smooth curve $\beta(t)$ with cusps on both 
	ends. The number $a=-2$  lies inside
	the unstable region according to 
	Figure \ref{fig2b}  and
	$a=(0.4+0.4\iota)$ which lies outside the 
	curve, is unstable as well (Figure \ref{fig2c}). 
	The system is unstable for any value of parameter $a$.
\end{Ex}
\begin{Ex}
	Now consider the region between $b=g_{1}(\alpha)$ and $b=g_{2}(\alpha)$. 
	Consider $b=2$ and $\alpha=0.5$. Here, 	the stability curve is divided into two parts (see Figure \ref{fig3a}). Consider a system (\ref{eqn2}) with $a=-1.2$.	It is inside the stable region and we get a stable trajectory as seen in Figure \ref{fig3b} 		For $a=0.2$, which lies in the unstable region 	we get an unstable solution (Figure 	\ref{fig3c}). The solutions for the values of $a$ lying outside the stability curve were checked and found unstable.
\end{Ex}
\begin{Ex}
	On the boundary $b=g_{2}(\alpha)$ which	is described by $b=1$ for any $\alpha$. Consider the stability curve for $\alpha=0.4$ and $b=1$(see Figure \ref{fig4a}). It generates a single region. The quantity $a=0.5$ inside the curve leads to a stable solution.		For $a=-0.7\iota$ outside the curve, an unstable solution is observed.
\end{Ex}
\begin{Ex}
	For $b$ in the region between $b=g_{2}(\alpha)$ and $b=g_{3}(\alpha)$, a single stable region (Figure \ref{fig5a}) is observed. For $\alpha=0.3$ and $b=0.1$, the value $a=0.3$ is  inside this region and gives a stable solution (Figure \ref{fig5b}). On the other hand, $a=-0.4$ lies outside the stable region and leads to an unstable solution (Figure \ref{fig5c}).
\end{Ex}
\begin{Ex}
	Examples with values of $b$ on the boundary curves give 
	limiting cases. For example $b=g_2(\alpha)=1$ is the limiting
	value of $b$ at which a single stable region is observed and
	for in the next region, we have one stable and one
	unstable region inside the stability curve.
	( See Example 3.3, 3.2). Similarly 
	for $b=g_3(\alpha)=-1$ the values for $t=0$ and $t=\pi$
	touch each on the real axis and we have 
	two stable regions. For  
	$b<g_3(\alpha)=-1$, we have 3 regions one of
	which is unstable.
	Consider the stability curve for $\alpha=0.7$ and $b=-1$ (see Figure \ref{fig6a}). Stable trajectories are observed for $a=1.1+0.5\iota$ and $a=1.3-0.8\iota$ (Figure \ref{fig6b} and \ref{fig6c} respectively). We checked the stability of the parameter values lying outside the curve and verified that they are unstable. 
\end{Ex}
\begin{Ex}
	Now we consider region lying  between $b=g_{3}(\alpha)$ and $b=g_{4}(\alpha)$. The stability curve for $\alpha=0.2$ and $b=-3$ is sketched as shown in Figure \ref{fig7a}. The curve has three regions, two of which look identical due to symmetry with respect to the real axis. The system with $a=(4+1.3\iota)$ has a stable solution (see Figure \ref{fig7b}). (Similarly, $a=(4-1.3\iota)$  leads to a stable trajectory.) Thus these two regions have stable solutions. On the other hand, for $a$ in the third region (say, $a=2$), we observe unstable trajectories (see Figure \ref{fig7c}). It was checked that	the parameter values  $a$ lying outside the stability curve lead to unstable solutions.
\end{Ex}
\begin{Ex}
	The last region lies below the curve $b=g_{4}(\alpha)$.	Consider $\alpha=0.5$ and $b=-2.2$, the curve $\beta$ (see Figure \ref{fig8a}) encloses  a single region. There are no stable values of $a$ inside or outside the region. Consider $a=1.5$ and $a=2.8$. Both values	lead to unstable trajectories (see Figure \ref{fig8b} and \ref{fig8c}, respectively). 
\end{Ex}

\section{Stability result for $\tau=2$}
\begin{figure}
	\includegraphics[width=\columnwidth]{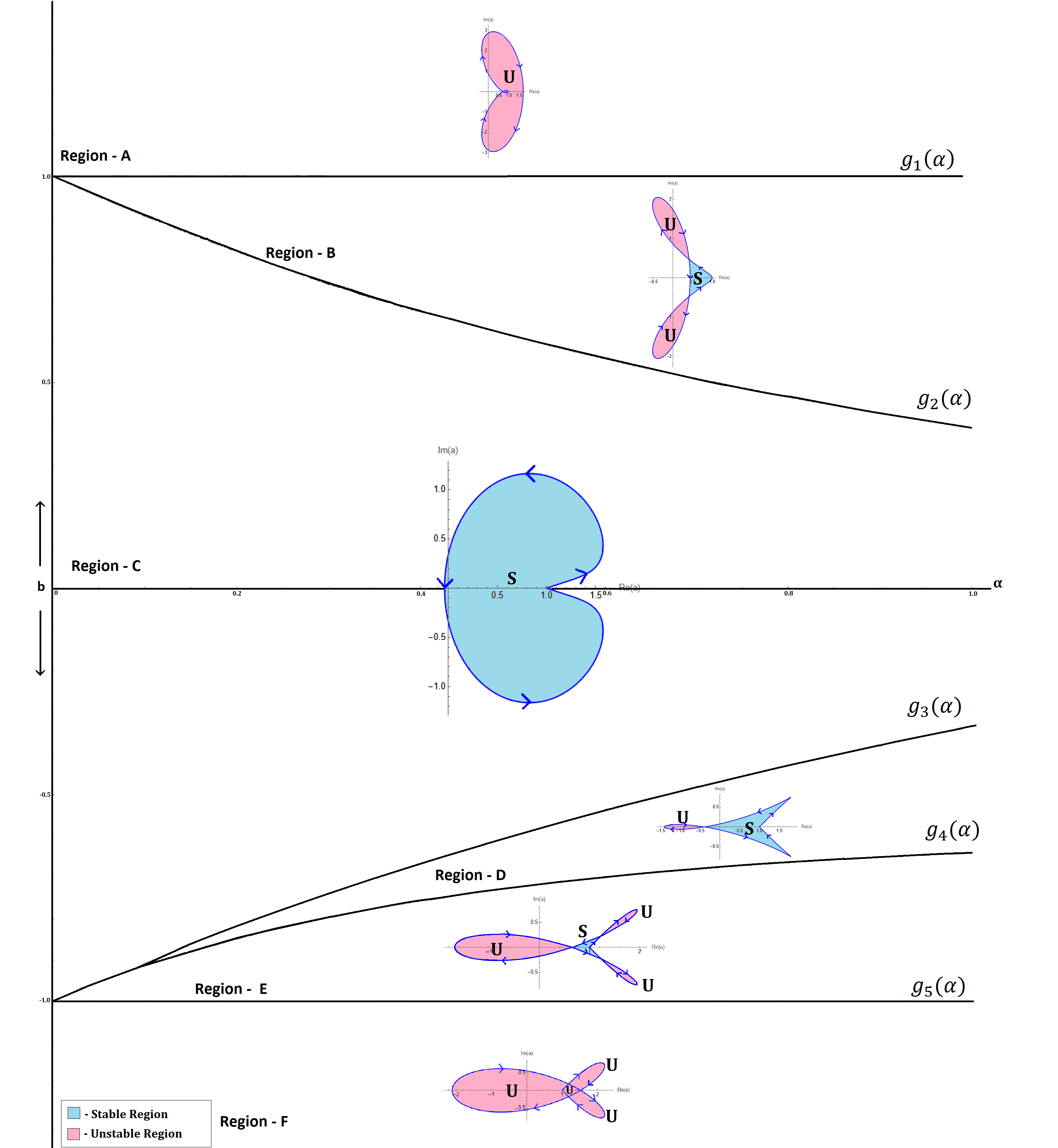}
	\centering
	\caption{bifurcation regions in $b \alpha$-plane for $\tau=2$ with representative stability diagrams}
	\label{fig9}
\end{figure}
\begin{The}
	Now consider the delayed fractional order equation (\ref{eqn2})	with $\tau=2$. The dynamics is richer and leads to five boundary curves that demarcate qualitatively different behaviors. Let $b=g_{j}(\alpha)$, $j=1,2,\cdots,5$ be the bifurcation		curves in the $b \alpha$-plane (cf. Figure \ref{fig9}) which are the branches of the implicit curve.
	\begin{equation}
		\begin{aligned}
			g(b,\alpha)=& \cos\left(\frac{1}{2}\left(\pi\alpha-(-3+\alpha)K\right)\right)
			-b(1+\alpha)\cos\left(\frac{1}{2}\left(\pi\alpha-(3+\alpha)K\right)\right)
			+b\sin\left(\frac{1}{2}(1+\alpha)\left(\pi-K\right)\right)\\
			& +(\alpha-1)\sin\left(\frac{1}{2}\left(\pi(1+\alpha)-(-1+\alpha)K\right)\right),
		\end{aligned} \label{cusp2}
	\end{equation}
	where, $K=\arccos\left(\frac{-2+2\alpha-\alpha^{2}+b^{2}(2+2\alpha+\alpha^{2})}{2(-1+\alpha+b^{2}(1+\alpha))}\right)$. \\
	For any $\alpha\in(0,1)$, $-1=g_{5}(\alpha)<g_{4}(\alpha)<g_{3}(\alpha)<g_{2}(\alpha)<g_{1}(\alpha)=1$.
	%The system has stable solutions for the boundaries mentioned below.
\end{The}
The stability results for the system  (\ref{eqn2}) with $\tau=2$ are as follows:\\
1. If $b \in (-\infty,g_{5}(\alpha)) \cup (g_{1}(\alpha),\infty)$,
then the system is unstable.\\
2. If $g_{5}(\alpha)<b<g_{4}(\alpha)$, then the boundary
curve $\beta(t)$ generates four regions; of which the outer three are unstable
and the central one is stable.\\
3. If $g_{4}(\alpha)<b<g_{3}(\alpha)$, then the number of
unstable regions reduces to one while the number of
the stable region remains the same, i.e. one. The two unstable regions described in Case 2 merge with the stable region and hence the size of the stable region, in this case, is larger than that of in Case 2.\\
4. If $g_{3}(\alpha)<b<g_{2}(\alpha)$, then the region
produced by the boundary curve is a single closed curve
which is stable.\\
5. If $g_{2}(\alpha)<b<g_{1}(\alpha)$, then the boundary
curve produces three regions, one stable and two unstable.\\
Now we consider the boundary cases:\\
1. If $b=g_{5}(\alpha)=-1$, then the curve $\beta(t)$ has
three regions and all the regions are unstable.\\
2. If $b=g_{4}(\alpha)$, then the curve $\beta(t)$ has
two regions, a stable region with two cusps on the right
side and an unstable region.\\
3. If $b=g_{3}(\alpha)$, then the curve $\beta(t)$ has
a single stable region with a cusp on the left side.\\
4. If $b=g_{2}(\alpha)$, then the curve $\beta(t)$ has
a single region that is stable and has two cusps.\\
5. If $b=g_{1}(\alpha)=1$, then the curve $\beta(t)$ has
two regions and both regions are unstable.\\

\textbf{Proof:}\\
As discussed in the proof of Theorem 3.1, as the parameter values change the simple closed curve $\beta(t)=(x(t), y(t))$ possesses multiple points and generates different bounded regions, in this case also. Among these, the regions with the anticlockwise oriented boundary are the stable ones. The multiple points are generated through the cusps.

To find the conditions for the $\frac{0}{0}$ form of $\lim_{t\rightarrow t_{0}}\beta'(t)$, we solve $x'(t)=0$ and $y'(t)=0$ simultaneously.
Squaring and adding these equations, we get

$$t=\arccos\left(\frac{2(\alpha-1)-\alpha^{2}+b^{2}(\alpha^{2}+2\alpha+2)}{2(b^{2}(1+\alpha)+\alpha-1)}\right)$$.

Substituting this value in the equation $x'(t)=0$, we get the expression for the existence of cusps in $\beta(t)$ as equation (\ref{cusp2}).
This implicit curve $g(b,\alpha)$ has 5 branches namely $b=g_{j}(\alpha)$, $j=1,2,\cdots,5$ as shown in Figure \ref{fig9}.
These bifurcation curves produce 6 different regions labeled as A, B, C, D, E, and F of $b- \alpha$ plane (see Figure \ref{fig9}). We list our observations below:\\

If the point $(b,\alpha)$ belongs to the region A {\it{i.e.}}, if $b>g_{1}(\alpha)$, the boundary curve $\beta(t)$ is a simple closed curve with a clockwise orientation. Thus, the regions inside as well as outside $\beta(t)$ are unstable.

In the region B, we have $g_{2}(\alpha)<b<g_{1}(\alpha)$. Here, $\beta(t)$ becomes a multiple curve and has two unstable and one stable bounded region. As $b$ decreases in B, the range of stable region goes on increasing.

If $g_{3}(\alpha)<b<g_{2}(\alpha)$ (i.e. region C), $\beta(t)$ becomes
a simple closed curve with anticlockwise orientation. Thus the inside part of $\beta(t)$ is the stable region and the outside part is unstable.
As the value of $b$ goes on decreasing, the stable region gets stretched horizontally and  a cusp is formed for $b=g_{3}(\alpha)$.

As $b$ is further reduced {\it{i.e.}}, when $g_{4}(\alpha)<b<g_{3}(\alpha)$ (i.e. region D),	the cusp on the boundary curve $\beta(t)$ evolves into a loop-shaped structure
of clockwise orientation  making $\beta$ a multiple curve with one stable
and one unstable region. As the value of $b$ moves closer towards $g_{4}(\alpha)$,
the area of unstable regions increases. At the bifurcation curve $b=g_{4}(\alpha)$,
two cusps are observed on the boundary curve $\beta$ on the right side.

In region E ($g_{5}(\alpha)<b<g_{4}(\alpha)$),  both of the cusps on the right of $\beta$ evolve into two unstable regions along with the previous
unstable region on the left side. Thus, the boundary curve $\beta(t)$ has 3 unstable regions and one stable region. Reducing the value of $b$ further results in the reduction of the size of the stable region.

When $b$ reaches the bifurcation curve $b=g_{5}(\alpha)$, the stable region disappears completely, leaving the boundary curve $\beta(t)$ with 3 unstable regions. These regions are unstable as a result of their clockwise orientation.

In region F, {\it{i.e.}} for $b<g_{5}(\alpha)$ the boundary curve has four bounded unstable regions with a clockwise orientation. Thus, below the bifurcation curve $g_{5}(\alpha)$, the system	becomes completely unstable for all values of $a$. 

As in the case of $\tau=1$, we have checked all the above results numerically.

\begin{figure*}[h]
	\centering
	\subfloat[]{%
		\includegraphics[scale=0.5]{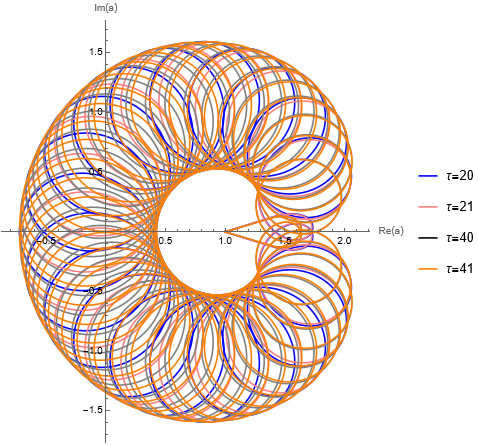}
	}
	\subfloat[]{%
		\includegraphics[scale=0.56]{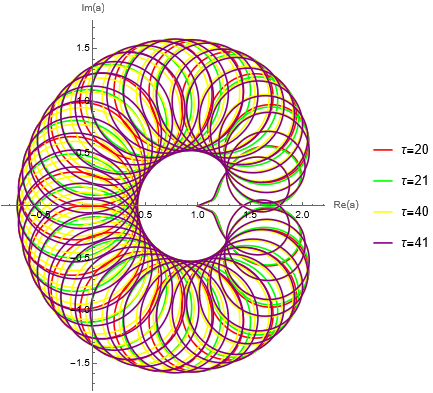}
	}
	\caption{Boundary curves are shown for the system (\ref{eqn 7}) with $\alpha=0.2$, $\tau=20, 21, 40, 41$, (a) $b=-0.5$ and (b) $b=0.5$.}
	\label{figd}
\end{figure*}

\section{Asymptotic limit of large delay}
We observe that we can stabilize a larger range of $a$ values with $\tau=1$ compared to $\tau=2$. The reason is that there are more routes to instability for larger $\tau$. We now consider the system (\ref{eqn 8}) for higher $\tau$ values. Let us take $\tau=20,21,40 $ and $41$ for $\alpha=0.2$. The boundary curves $\beta$ have been plotted for the values of $b=-0.5$ and $b=0.5$ in the figures \ref{figd} (a) and (b) respectively. In both figures, the stable region is in the center and is surrounded by the spring-like structure. This spring-like structure has a clockwise orientation and thus it is unstable. On the real axis, the upper bound is one. The lower bound is indicated intersection of the stability curve with the real axis.
The lower bound of real stable $a$ is at $\beta(\pi)$ when $\tau$ is even and $b$ is positive or $\tau$ is odd and $b$ is negative. However, as we increase $\tau$, the lower bound approaches $\beta(\pi)$.  For a given $\tau$, the stable region  depends slightly on whether $\tau$ is even or odd and the sign of $b$. This dependence reduces as we increase $\tau$. The upper bound is fixed at $t=0$ for the boundary curve $\beta(t)$, the lower bound approaches $\beta(\pi)$. We find that for large $\tau$, the stable range of $a$ does not depend on the sign of $b$ and is given by $a=1+2^{\alpha}(\vert b \vert -1)$.

The $b-a$ curves plotted in the figures \ref{figb1} and \ref{figb} can be used to determine the values of feedback coefficient $b$ to get the maximum stable range. As $\tau$ increases the $b-a$ curve reaches a limiting region. As mentioned above, the stable $b-a$ region  is given by $a=1$ and $a=1+2^{\alpha}(\vert b \vert -1)$ as $\tau\rightarrow\infty$. 

In a control system, the feedback term can be selected according to the requirements of the systems. We see that the maximum range of asymptotically stable fixed points can be obtained in the case of $\tau=1$. So, one can use $\tau=1$ to maximize the range of asymptotically stable fixed points of the system. Large delay can be used if the system requires enhanced chaos.

\section{Efficacy of control for large negative multipliers}
\begin{figure*}[h]
	\centering
	\includegraphics[scale=0.6]{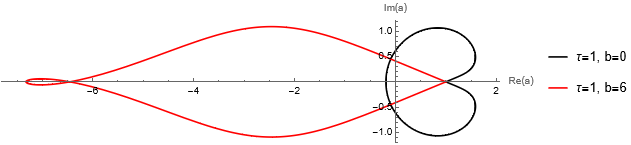}
	\caption{Stability regions for system (\ref{eqn2}) $\alpha=0.25$.}
	\label{fige}
\end{figure*}

\begin{figure*}[h]
	\centering
	\subfloat[Unstable solution for $b=0$.]{%
		\includegraphics[scale=0.4]{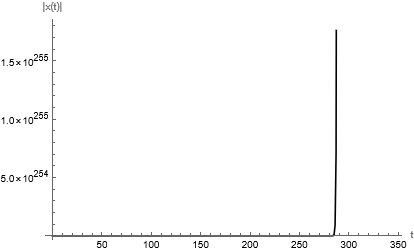}
	}\hspace{0.3cm}
	\subfloat[Stable solution for $b=6$.]{%
		\includegraphics[scale=0.4]{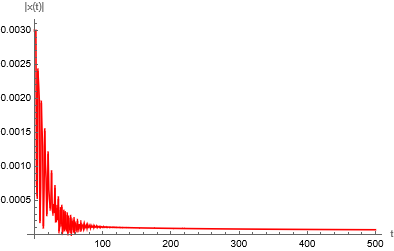}
	}
	\caption{Stability analysis for system (\ref{eqn2}) with $\alpha=0.25$, $\tau=1$ and $a=-7$.}
	\label{figf}
\end{figure*}

The range of stable real value of $a$ ranges from $1-2^{\alpha}$ to 1 for system (\ref{eqn1}). The real range for integer order difference equation is between -1 and 1. There is a reduction in the real range for fractional order maps. Now $a=-7$ is outside stable range for any $0<\alpha\le1$. The stability of system (\ref{eqn2}) depends on the feedback coefficient ‘$b$’ and the delay ‘$\tau$’. For $b=0$ , we have system (\ref{eqn1}) without feedback. Figure \ref{fige} shows stability curves for system (\ref{eqn2}) with $\alpha=0.25,\tau=1, b=0$ (in black) and $b=6$ (in red). The point $a=-7$ lies outside the black stability curve. It is enclosed within small stable region of the red curve. We analyze the stability of this point in both cases. We iterate the systems for $t=500$ time-steps. For $b=0$ case, the trajectory diverges and the solution is unstable (see Figure \ref{figf}(a)). For $b=6$ case, the trajectory asymptotically converges to zero and the solution is stable(see Figure \ref{figf}(b)). Thus large negative multipliers can be stabilized with feedback for $\tau=1$.

\begin{figure*}[h]
	\centering
	\includegraphics[scale=0.55]{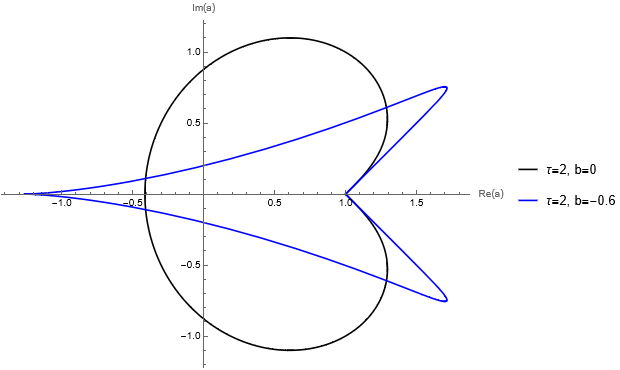}
	\caption{Stability regions for system (\ref{eqn2}) $\alpha=0.5$.}
	\label{figg}
\end{figure*}

\begin{figure*}[h]
	\centering
	\subfloat[Unstable solution for $b=0$.]{%
		\includegraphics[scale=0.43]{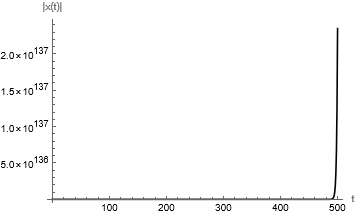}
	}\hspace{0.3cm}
	\subfloat[Stable solution for $b=-0.6$.]{%
		\includegraphics[scale=0.43]{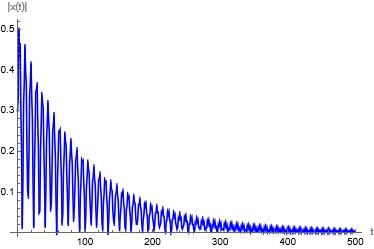}
	}
	\caption{Stability analysis for system (\ref{eqn2}) with $\alpha=0.5$, $\tau=2$ and $a=-1.1$.}
	\label{figh}
\end{figure*}

For the system (\ref{eqn2}) with $\tau=2$, the effect is less dramatic. Figure \ref{figg} shows the stability curves for systems (\ref{eqn2}) with $\alpha=0.5,\tau=2, b=0$ (in black) and $b=-0.6$ (in blue). In the first case $b=0$, the point $a=-1.1$ lies outside the black stability curve and has an unstable solution (see Figure \ref{figh}(a)). In second case $b=-0.6$, the point $a=-1.1$ lies inside the blue stability curve and has a stable solution as seen in Figure \ref{figh}(b).

Thus, it is possible to stabilize a larger range of parameter values for a fractional order system by adding a feedback term with delay (system (\ref{eqn2})). The case of $\tau =1$ gives the maximum range of stability with appropriate $b$ and is recommended for control.

\section{Nonlinear Maps}

Consider a fractional order map with a feedback control
\begin{equation}
	x(t)=x_0 + bx(t-\tau)+\frac{1}{\Gamma(\alpha)}\sum_{j=1}^{t}\frac{\Gamma(t-j+\alpha)}{\Gamma(t-j+1)}[f(x(j-1))-x(j-1)]. \label{eqn 8}
\end{equation}
\begin{Def}
	A steady-state solution of (\ref{eqn 8}) is called an equilibrium point. Thus, if $x_*$ is an equilibrium point of (\ref{eqn 8}) then $x_0=x_*$ implies $x(t)=x_*$ for all $t=1,2,\cdots$.
\end{Def}
Note that, the only equilibrium point of (\ref{eqn 8}) is $x_*=0$ when we have $f(0)=0$. Therefore, we consider the systems (\ref{eqn 8}) with equilibrium point $x_{*}=0$.

In a neighborhood of the point $x_*=0$, we have $f(x)\approx f(x_*)+(x-x_*)f'(x_*)=a x$, where $a=f'(0)$. Thus, the local stability properties of the nonlinear system (\ref{eqn 8}) at $x_*=0$ are the same as those of the linear system (\ref{eqn2}).

We consider only real values of $a$ and $b$, in this case. Therefore, the stable region given by the curve $\beta(t)$ defined by (\ref{eqn 7}) is bounded by the curves $a=1$, $a=-2^{\alpha}\left(1-b(-1)^{\tau}\right)$ and the parametric curve $a(t)=2^{\alpha}\left(\sin\left(\frac{t}{2}\right)\right)^{\alpha}\left(\cos\left(\frac{\alpha\pi}{2}+t\left(1-\frac{\alpha}{2}\right)\right)-\sin\left(\frac{\alpha\pi}{2}+t\left(1-\frac{\alpha}{2}\right)\right)\cot\left(\frac{\alpha\pi}{2}+t\left(1-\tau-\frac{\alpha}{2}\right) \right)\right)+1,$\\
$b(t)=\frac{\sin\left(\frac{\alpha\pi}{2}+t\left(1-\frac{\alpha}{2}\right)\right)}{\sin\left(\frac{\alpha\pi}{2}+t\left(1-\tau-\frac{\alpha}{2}\right)\right)},$ $t\in[0,2\pi]$ in the $b-a$ plane. We denote these curves as $b-a$ curves.

Thus, the system  (\ref{eqn 8}) can be controlled if the point $(b,a)$ with $a=f'(0)$ lies inside the region bounded by these curves.
We consider the nonlinear systems  (\ref{eqn 8}) of fractional order $\alpha$ and delay $\tau$.
The fractional logistic map in this model is defined by the equation  (\ref{eqn 8}) with $f(x)=\lambda x(1-x)$, where $\lambda$ is a parameter and $a=\lambda(1-2x_*)=\lambda$.\\
Similarly, the fractional cubic map is  (\ref{eqn 8}) with $f(x)=\beta x^{3}+(1-\beta)x$, where $a=3\beta x_*^{2}+(1-\beta)=1-\beta$.

\begin{figure*}[h]
	\centering
	\subfloat[Logistic map with $\alpha=0.25$.]{%
		\includegraphics[scale=0.2]{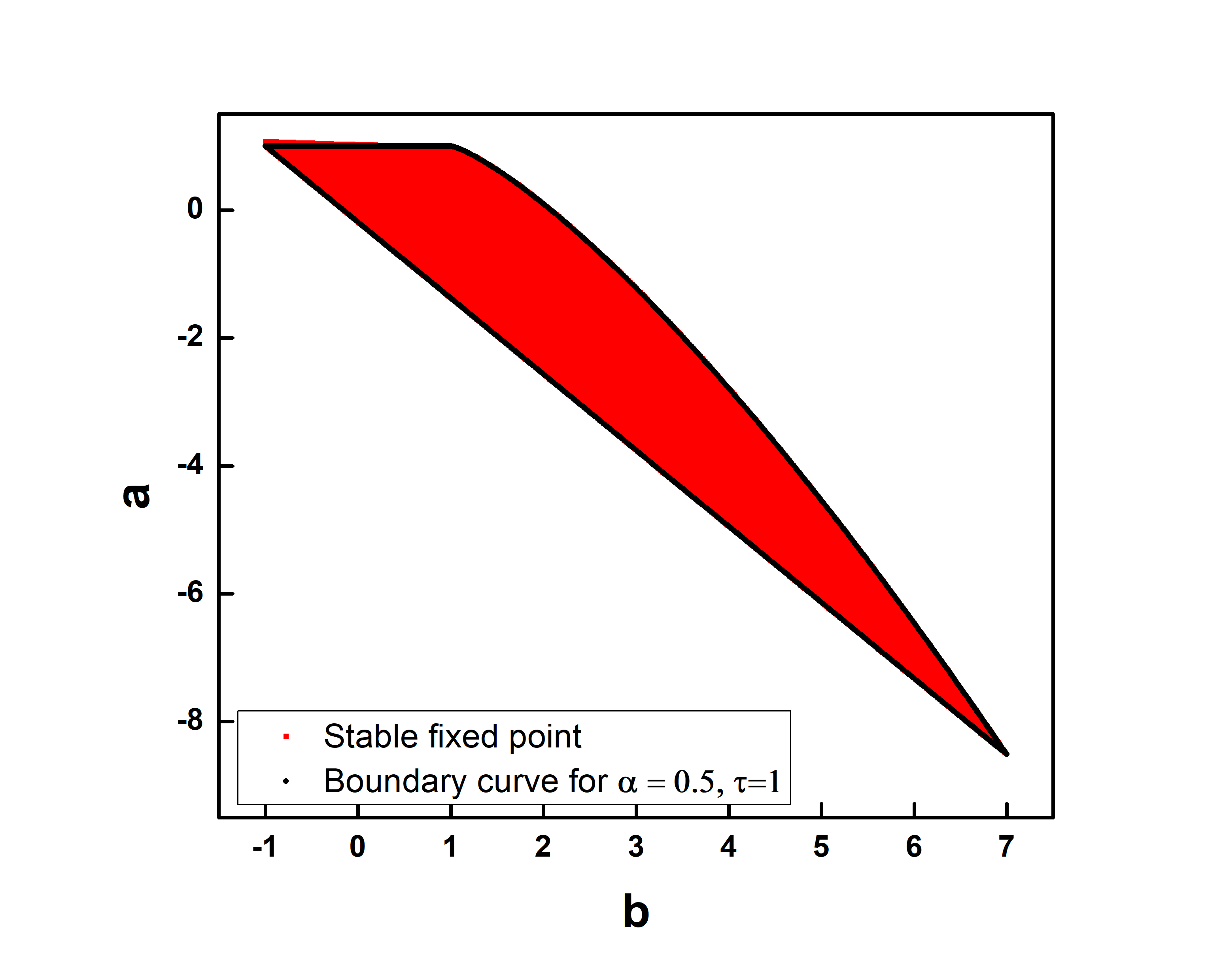}
	}
	\subfloat[Logistic map with $\alpha=0.5$.]{%
		\includegraphics[scale=0.2]{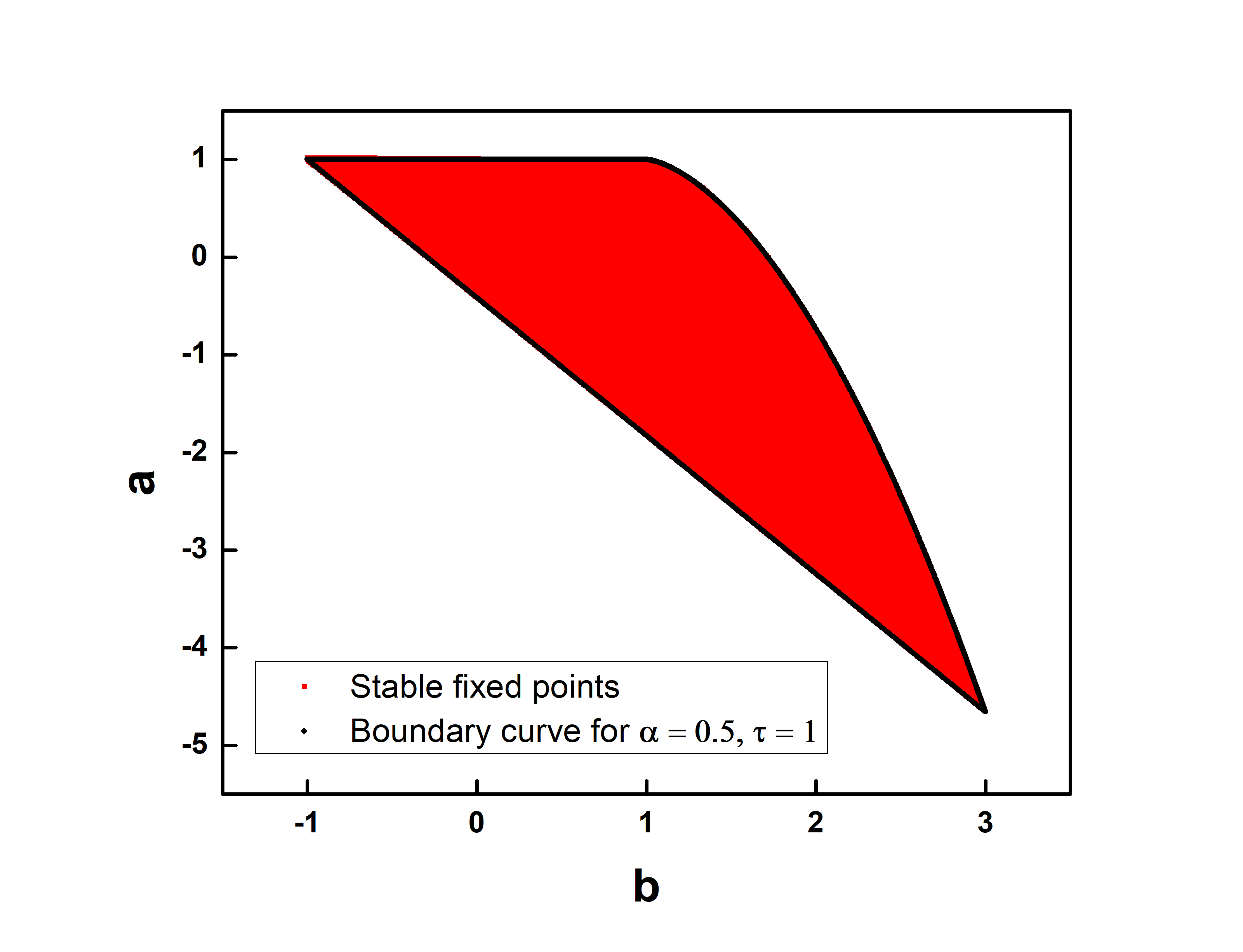}
	}
	\subfloat[Logistic map with $\alpha=0.75$.]{%
		\includegraphics[scale=0.2]{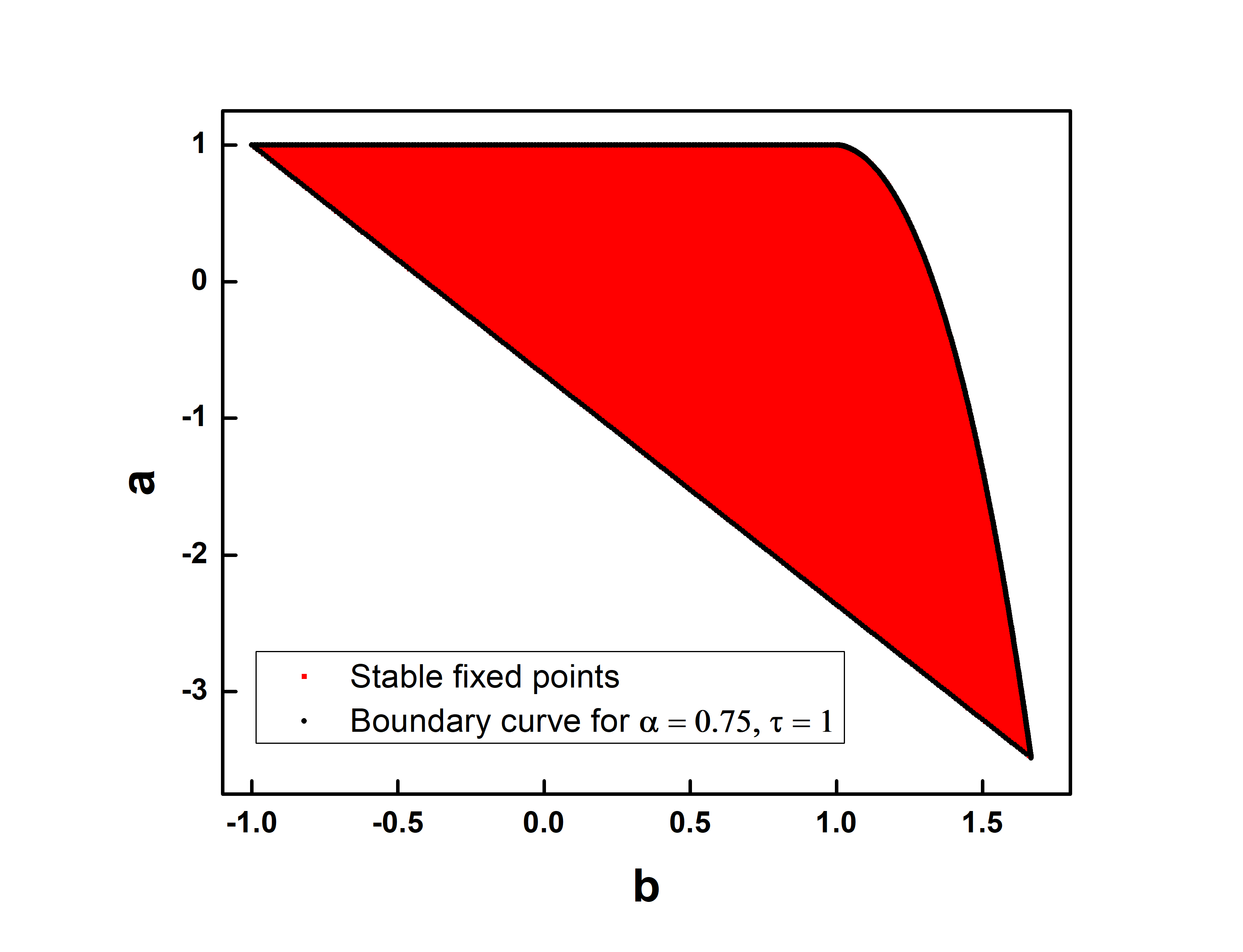}
	}\\
	\subfloat[Cubic map with $\alpha=0.25$.]{%
		\includegraphics[scale=0.2]{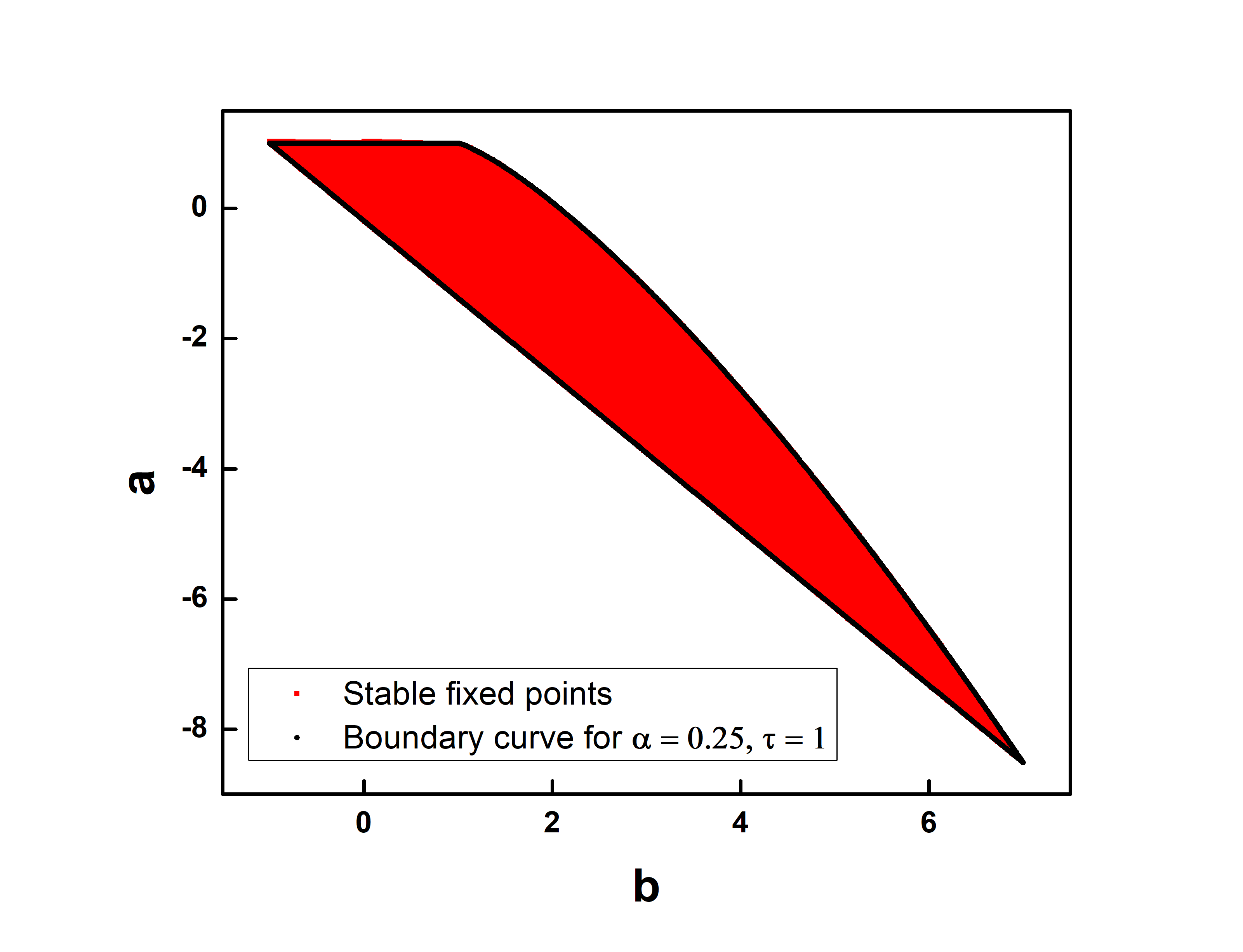}
	}
	\subfloat[Cubic map with $\alpha=0.5$.]{%
		\includegraphics[scale=0.2]{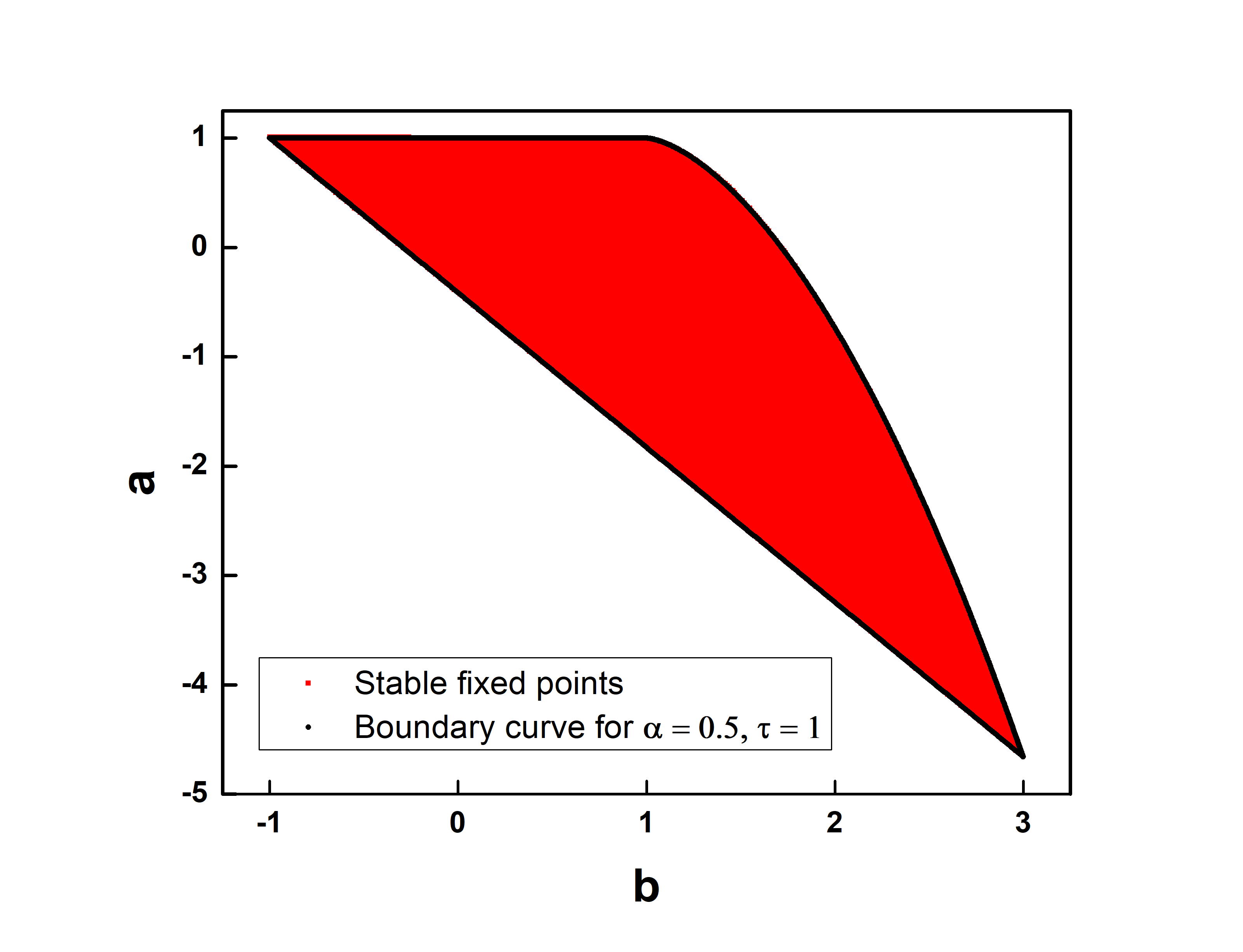}
	}
	\subfloat[Cubic map with $\alpha=0.75$.]{%
		\includegraphics[scale=0.2]{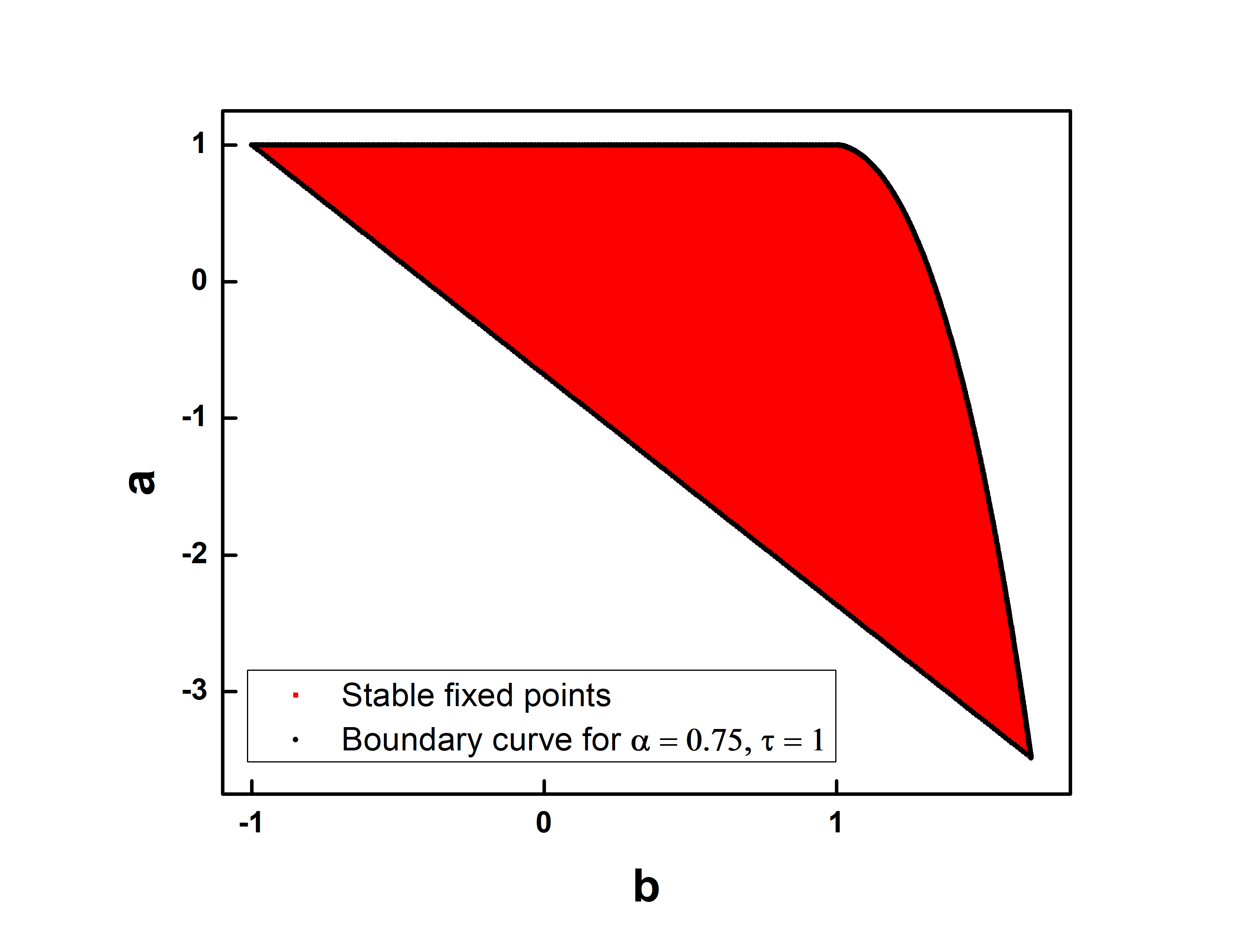}
	}
	\caption{Stable region for fractional systems with $\tau=1$. The stable fixed points of fractional systems lie within the region enclosed by the $b-a$ curve. We discard $T$ time-steps and the convergence is within $\delta$. For $\alpha=0.25$, $T=6\times 10^{4}$ and $\delta=10^{-9}$ for logistic map and $T=4\times 10^{4}$ and $\delta-5\time 10^{-8}$ for	cubic map. In all other cases $T=2\times 10^4$ and $\delta=10^{-5}$.}
	\label{figb1}
\end{figure*}

Figure \ref{figb1} shows the stability regions of all these systems for the $\tau=1$ for various values of $\alpha$.

We iterate these maps for $T$ time-steps where $T$ is large and the equilibrium or asymptotic equilibrium point is assumed to be stable if the convergence is obtained within $\delta$. We know the boundaries of the stable fixed point in the $b-a$-plane a for linear system. For $a=f'(0)$, we observe that the $(b,a)$ values for stable fixed point lie in the analytically obtained bounds for both $\tau=1$ and $\tau=2$. For $\tau=1$, we have plotted $(b,a)$ values in Figure \ref{figb1} for $\alpha=0.25,0.5,0.75$. It is observed for all the stable zero fixed points of the logistic, and cubic maps lie within the stable region defined by the $b-a$ curves.

We find that this $b-a$ region indeed encloses the stability region and gives precise bounds. Thus if the stable equilibrium point is known, the required strength of feedback for which the equilibrium point is stabilized can be found. The uncontrolled case corresponds to $b=0$ and if we want to stabilize (or destabilize) the equilibrium point an appropriate value of $b$ can be used and desired range of stability can be obtained.  We have studied $\alpha=0.25,0.5,0.75$ and $\tau=1$. For $\tau=1$, the stability region decreases in size with an increasing $\alpha$.

\begin{figure*}[h]
	\centering
	\subfloat[Logistic map.]{%
		\includegraphics[scale=0.23]{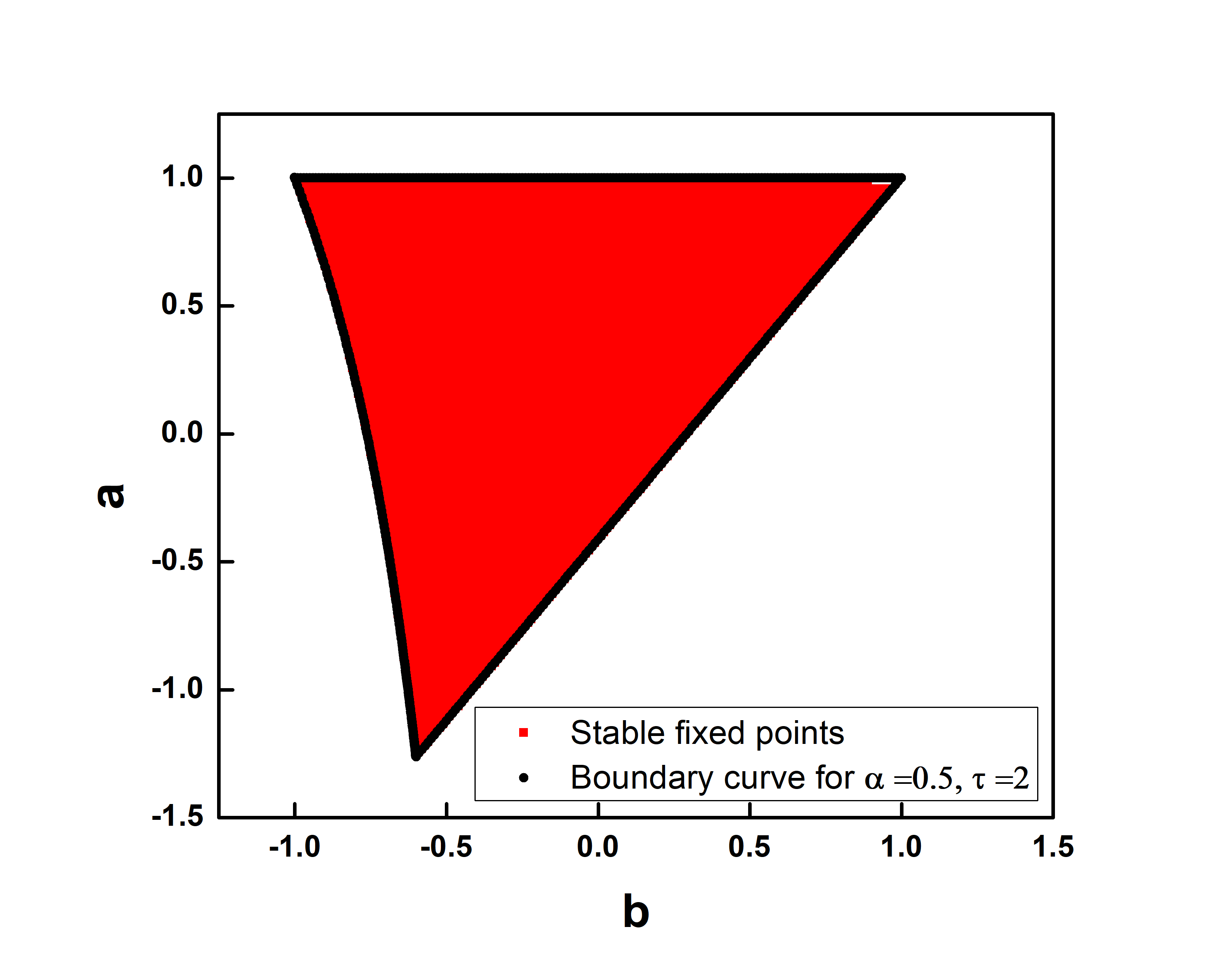}
	}
	\subfloat[Cubic map.]{%
		\includegraphics[scale=0.23]{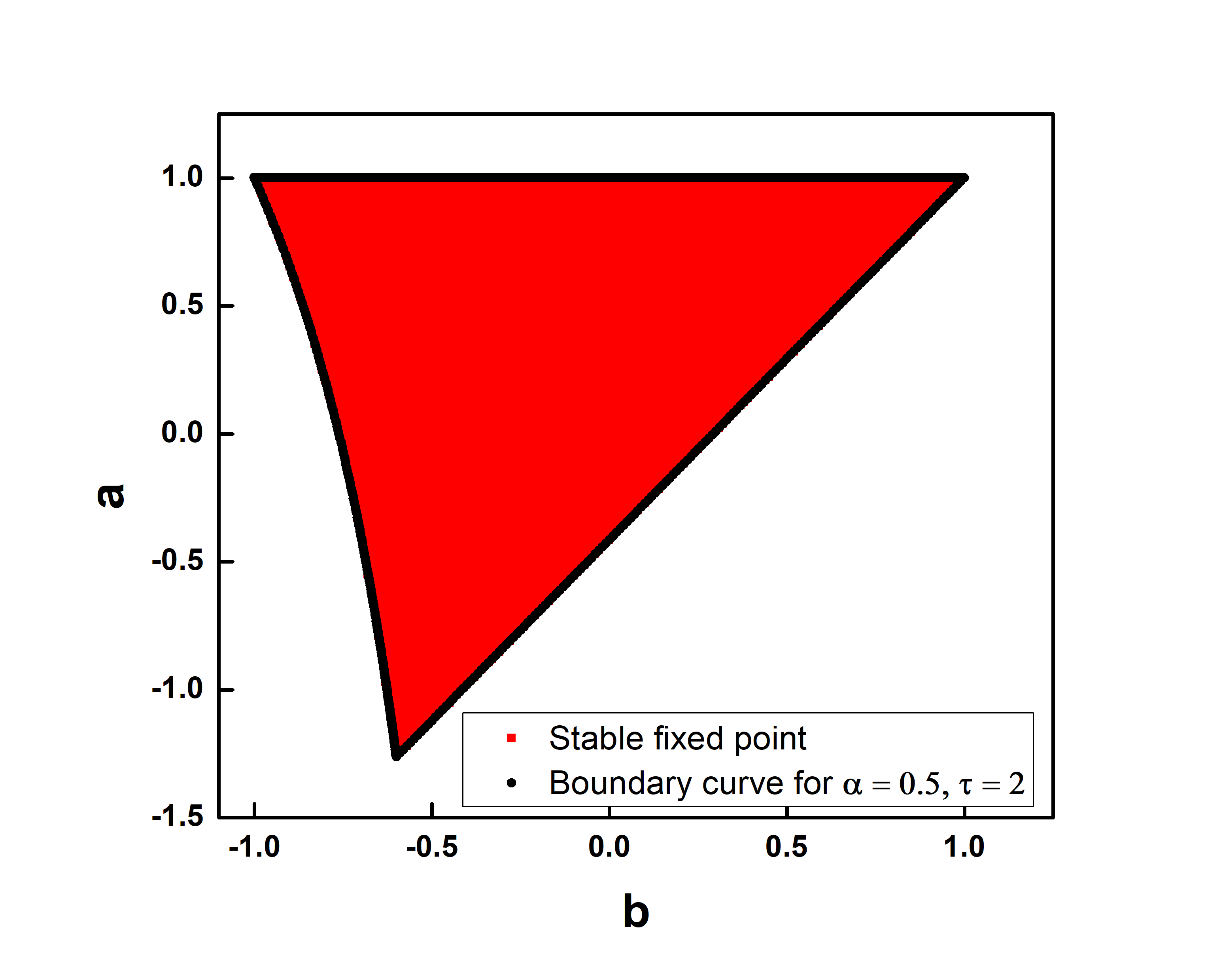}
	}
	\caption{Stable region for fractional systems with $\alpha=0.5$, $\tau=2$. The stable fixed points of the fractional systems lie within the region enclosed by the $b-a$ curve. We discard $T$ time-steps and the convergence is	within $\delta$.  For logistic map $T=10^5$ and $\delta=10^{-9}$ and for cubic map $T=5\times 10^{4}$ and $\delta=5\times 10^{-7}$.}
	\label{figb}
\end{figure*}

Figure \ref{figb} confirms that the stable fixed point of the fractional maps lies within the curve enclosed by the $\alpha$ dependent $b-a$ curves. Similarly, for $\tau=2$, stability regions can be obtained for various values of $\alpha$.	The stable regions for the zero fixed point for logistic, and cubic map with fractional order $\alpha=0.5$ and $\tau=2$  lie within the stable region derived for linear case (see Figure \ref{figb}).

\begin{figure*}[h]
	\includegraphics[scale=0.5]{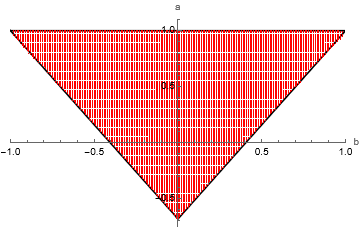}
	\centering
	\caption{Stability region in $b-a$ plane for the system (\ref{eqn 8}) with $\alpha=0.75$ and large $\tau$. The asymptotic stable fixed points of the fractional logistic map defined by (\ref{eqn 8}) with $\alpha=0.75$ and $\tau=100$ lie inside the generalized stability region.}
	\label{figd1}
\end{figure*}

Now we consider the case of large $\tau$. We consider the fractional logistic map defined by the equation (\ref{eqn 8}) with $f(x)=\lambda x(1-x)$, where $f'(0)=\lambda=a$. We iterate this map for  $T=6\times10^{4}$ and the equilibrium point is assumed to be asymptotically stable if convergence is obtained within $\delta=10^{-7}$. Here, we consider $\alpha=0.75$ and $\tau=100$. The triangular stability region encloses the asymptotically stable fixed points for the fractional logistic map. This is the generalized stable region for large $\tau$. From the Figure \ref{figd1} it is clear that  the range of stable region does not increase for any nonzero  value of $b$. It decreases with larger $\vert b\vert$. Thus, the large delay is not useful in stabilizing the fixed point  for the fractional system (\ref{eqn 8}) and
the best results are obtained with $\tau=1$. However, as mentioned above, larger $\tau$ can be used to destabilize the stable system. 

\begin{figure*}[h]
	\centering
	\subfloat[]{%
		\includegraphics[scale=0.25]{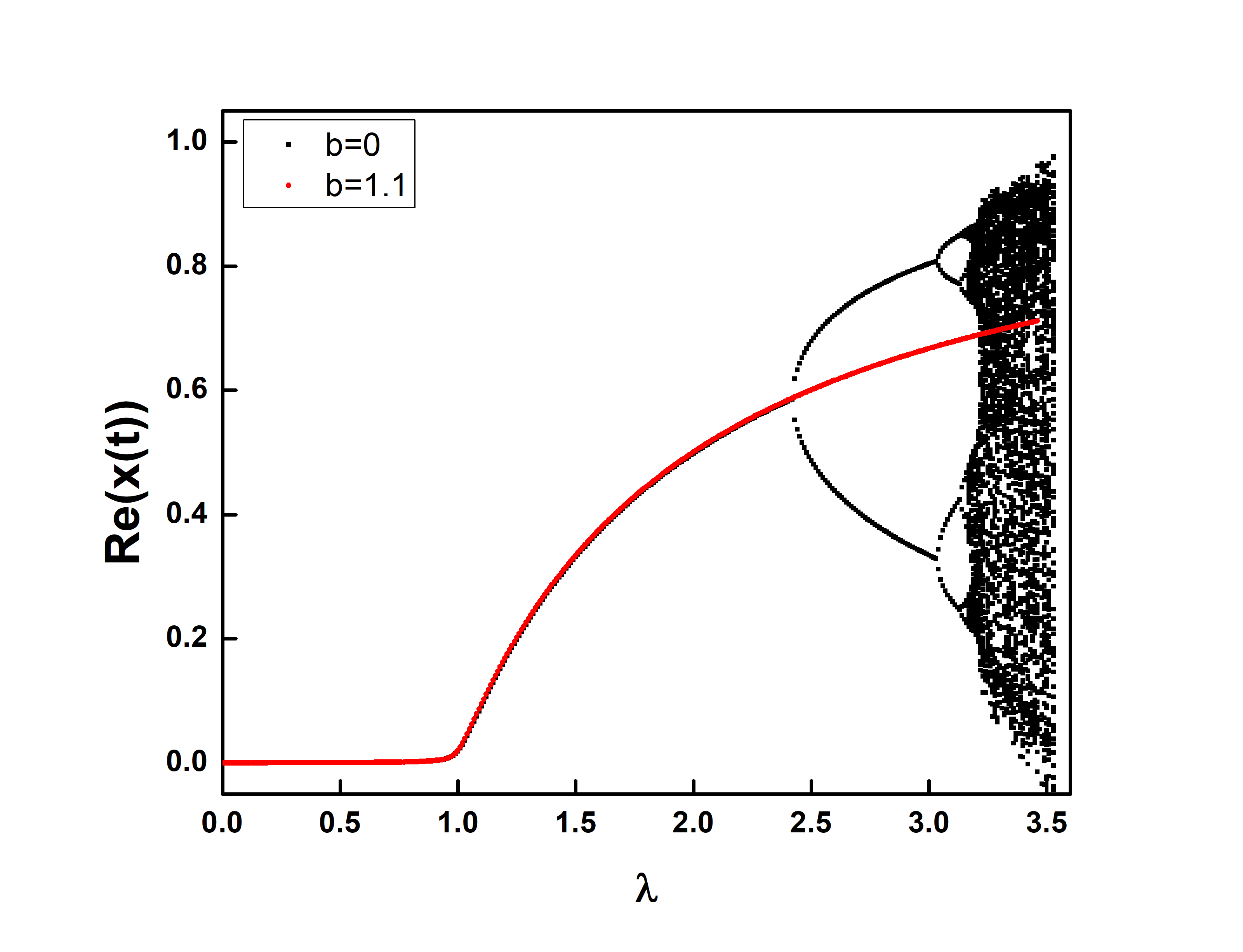}
	}
	\subfloat[]{%
		\includegraphics[scale=0.25]{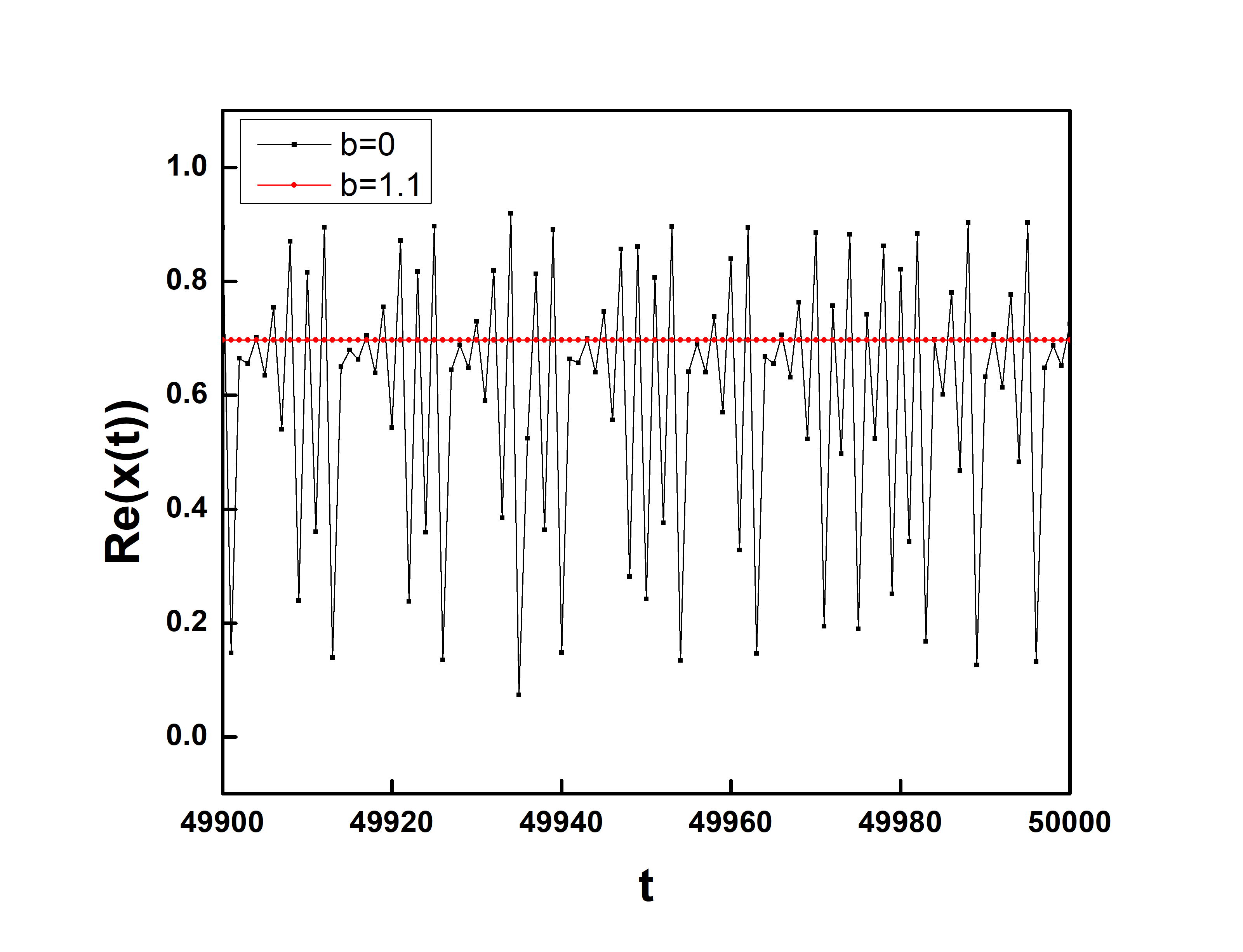}
	}\\
	\subfloat[]{%
		\includegraphics[scale=0.25]{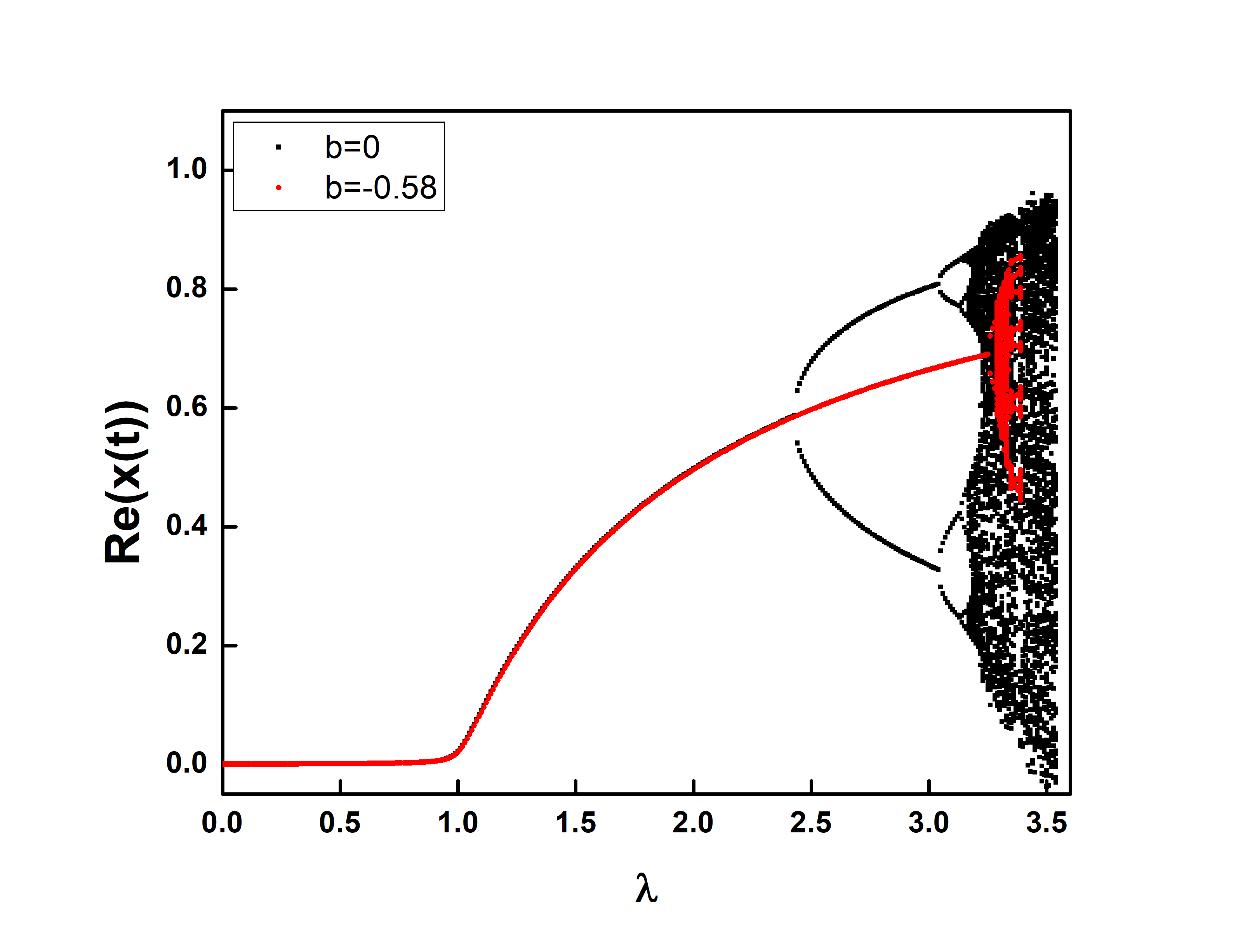}
	}
	\subfloat[]{%
		\includegraphics[scale=0.25]{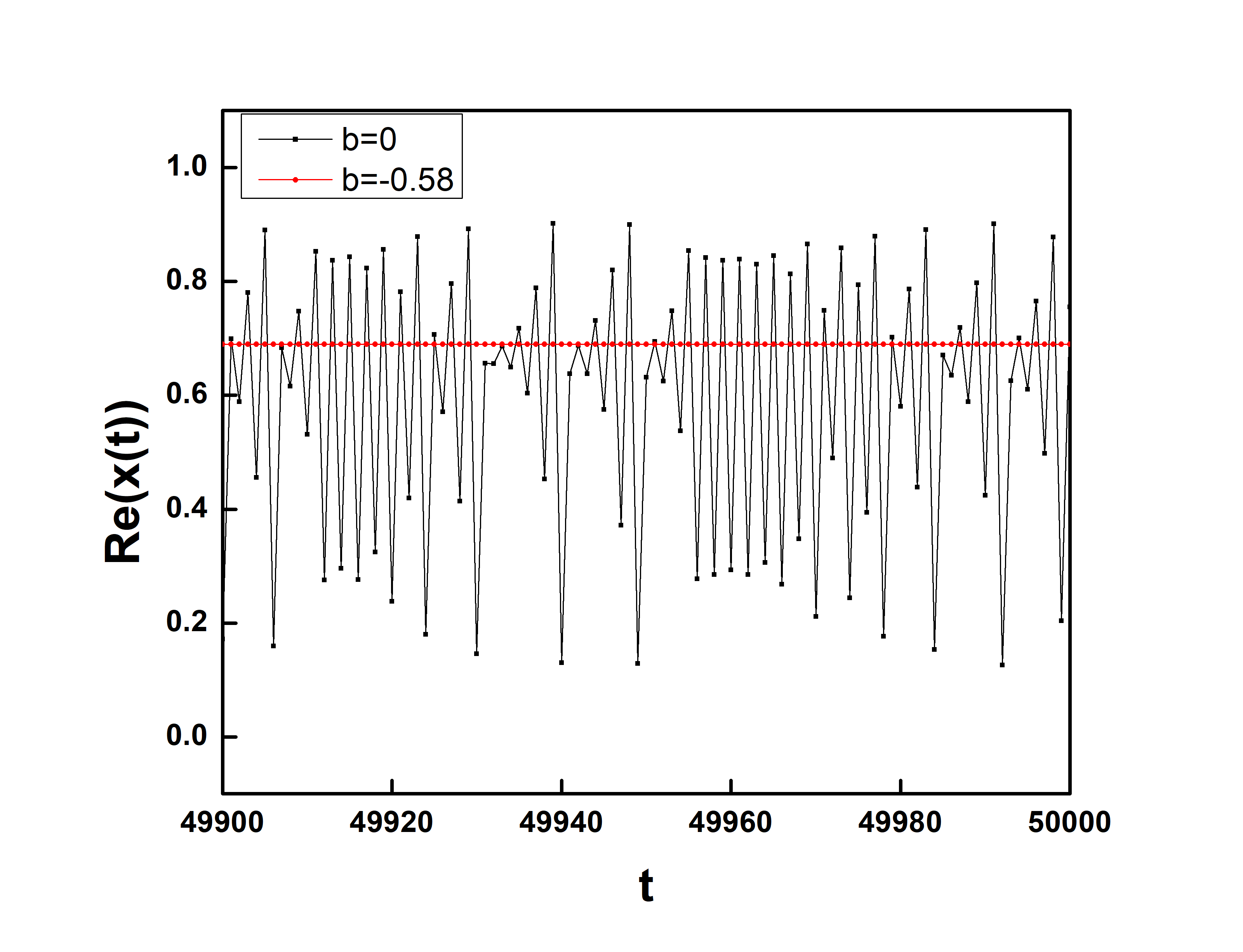}
	}
	\caption{Figure (a)  shows the bifurcation diagram for a logistic map with $\tau=1$ and $b=1.1, \alpha=0.5$ where only a stable fixed point is observed. The bifurcation diagram without control is also shown for reference. Figure (c) shows bifurcation diagram for $\tau=2$ and $b=-0.58,\alpha=0.5$. The bifurcation diagram without control is also shown for reference. The range over which a stable fixed point is observed in either case. Figure (b) shows the stabilization of the unstable chaotic state for $\tau=1$, $\alpha=0.5$, and $\lambda=3.3$ to the stabilized fixed point. Similarly, figure (d) shows the stabilization of the unstable chaotic state to stabilized fixed point for $\tau=2$, $\alpha=0.5$, and $\lambda=3.23$.}
	\label{figc}
\end{figure*}

The control term in this model is a delay term with coefficient `$b$'. In figure \ref{figc} (a) and (c), an increment in the range of nonzero stable fixed points can be seen. We observe that the proper selection of control parameter $b$ results in an extension of the range of the stable fixed points of the system. The maximum range of the nonzero stable fixed points can be found with the help of the $b-a$ curve for that particular $\alpha$ value the of fractional system. The range of the stability of nonzero fixed points is extended in presence of feedback. Preliminary numerical investigations indicate	that the above conditions are necessary but not sufficient for the	stability of nonzero fixed points. Further analytic investigations are needed to find the stability of nonzero fixed points. Figures \ref{figc} (b) and (d) show the time series of the system.

\section{Discussion and Conclusion}
Fractional order systems with delay have been investigated in many contexts.
In neural networks, numerous works have been done. In \cite{xu2021bifurcation}, stability conditions, and the existence of Hopf bifurcation for fractional order BAM neural network (FOBAMNN) with delay have been established. In \cite{huang2021bifurcations}, FOBAMNN with four delays has been turned into two delays and a correlation between stability and delay term has been studied. Investigation of fractional order neural network with multiple leakage delay shows that both fractional order and time delay are very important in controlling the transient behaviors of the FONN devised in \cite{huang2020bifurcations}. The stability of fractional order triangle multi-delayed neural networks has been studied in \cite{xu2022bifurcation}. A comparative study of integer order and fractional order delayed BAM neural network shows increased stability region \cite{xu2022comparative}. In \cite{xu2023new}, sufficient conditions for different delays for ensuring stability and generation of Hopf bifurcation have been demonstrated. Stability and bifurcation of an isovalent version of a fractional-order stage-structured predator–prey system have been investigated in \cite{xubifurcation}. Global asymptotic stabilization of fractional-order memristor-based neural networks (FMNNs) with time delay can be achieved by adjusting two groups of parameters as illustrated in \cite{jia2019global}. The analysis of fractional order Bloch equation with delay shows behaviors ranging from damped oscillations to oscillations with increasing amplitude for various values of delay \cite{bhalekar2011fractional}.  In this work, we carry out basic investigations in 
the context of fractional order diffference equations.

Control of chaos in dynamical systems is an important aspect of studies in the theory of dynamical systems from viewpoint of applications. In integer order differential and difference equations, it is a well-studied problem both experimentally and theoretically. One of the simplest control schemes is feedback and the Pyragas method is one such method with feedback delay control to stabilize the chaotic systems. In this work, we studied systems defined by fractional difference equations coupled with a delay term. The delay term acts as a control in this system. We give analytic conditions for the stability of the fixed point of these systems for the arbitrary delay. More detailed analysis is carried out for $\tau=1$ and $\tau=2$ followed by analysis for asymptotic limit. A detailed analysis is carried out for $\tau=1$ and $\tau=2$. We give all bifurcation curves. They are given by $g(b,\alpha)$. For real maps, we give the generalized stable regions for $\tau=1$, $\tau=2$ and in the asymptotic limit. Using $g(b,\alpha)$, we obtain the upper and lower bound of $b$ for particular $\alpha$ for which the stable regions exist. The stability region of a system depends on the fractional order $\alpha$, delay $\tau$, and the control parameter $b$. For nonlinear map $f(x)$, the zero  fixed point with slope '$a=f'(0)$' in the same $b-a$ range. The boundary curve (\ref{eqn 7}) may enclose stable and/or unstable region/s for multiple combinations of $\alpha$, $\tau$ and $b$. The criterion for finding out the stable and unstable region is also given using the orientation of the curve. Finally, we have studied nonlinear systems with fixed point $x_*=0$. Thus the results apply to a broad range of systems and this is a practical scheme for control. Our analysis can be used to stabilize/destabilize the fractional order difference system by introducing the appropriate delay term. The chaotic systems can be stabilized by selecting the delay $\tau=1$ whereas the stable system can be made unstable/chaotic by selecting a larger value for $\tau$.

\section{Acknowledgement}
PMG and DDJ thank DST-SERB for financial assistance (Ref. CRG/2020/003993).

\bibliographystyle{plain}
\bibliography{newref}

\end{document}